\documentclass[11pt]{article}
\usepackage{amssymb,amsmath,latexsym}
\usepackage{epsfig}

\newtheorem{thm}{Theorem}[section]
\newtheorem{cor}[thm]{Corollary}
\newtheorem{lemma}[thm]{Lemma}
\newtheorem{prop}[thm]{Proposition}

\numberwithin{equation}{section}

\def\nn{\nonumber}

\def\pf{{\medskip\noindent {\bf Proof. }}}
\def\qed{{\hfill $\Box$ \bigskip}}

\def\sA {{\cal A}}

  \def\sL {{\cal L}}

\def\bP {{\mathbb P}}  \def\bR {{\mathbb R}}

\def\R {{\mathbb R}}
\def\wt{\widetilde}
\def\wh{\widehat}

\def\E{{\mathbb E}}
\def\P{{\mathbb P}}

\def\bea{\begin{align*}}
\def\eea{\end{align*}}
\def\bee{\begin{equation}}
\def\eee{\end{equation}}

\def\eps{\varepsilon}

\def\wh{\widehat}

\makeatletter
\@addtoreset{equation}{section}

\makeatother

\begin{document}
\bibliographystyle{plain}

\title{\Large \bf
Dirichlet Heat Kernel Estimates for $\Delta ^{\alpha /2}+ \Delta
^{\beta /2}$}

\author{{\bf Zhen-Qing Chen}\thanks{Research partially supported
by NSF Grant DMS-0906743.}, \quad {\bf Panki Kim}\thanks{Research
supported by National Research Foundation of Korea Grant funded by
the Korean Government (2009-0087117).} \quad and \quad {\bf Renming
Song} }
\date{(October 16, 2009)}

\maketitle

\begin{abstract}
For $d\geq 1$ and $0<\beta<\alpha<2$, consider a family of pseudo
differential operators $\{\Delta^{\alpha} + a^\beta
\Delta^{\beta/2}; \ a \in [0, 1]\}$ on $\R^d$ that evolves
continuously from $\Delta^{\alpha/2}$ to $ \Delta^{\alpha/2}+
\Delta^{\beta/2}$. It gives arise to a family of L\'evy processes
\{$X^a, a\in [0, 1]\}$ on $\R^d$, where each $X^a$ is the sum of
independent a symmetric $\alpha$-stable process and a symmetric
$\beta$-stable process with weight $a$. For any $C^{1,1}$ open set
$D\subset \R^d$, we establish explicit sharp two-sided estimates
(uniform in  $a\in [0,1]$) for the transition density function of
the subprocess $X^{a, D}$ of $X^a$ killed upon leaving the open set
$D$. The infinitesimal generator of $X^{a, D}$ is the non-local
operator $\Delta^{\alpha} + a^\beta \Delta^{\beta/2}$ with zero
exterior condition on $D^c$. As  consequences of these sharp heat
kernel estimates, we obtain uniform sharp Green function estimates
for $X^{a, D}$ and uniform boundary Harnack principle for $X^a$ in
$D$ with explicit decay rate.
\end{abstract}

\bigskip
\noindent {\bf AMS 2000 Mathematics Subject Classification}: Primary
60J35, 47G20, 60J75; Secondary 47D07

\bigskip\noindent
{\bf Keywords and phrases}: fractional Laplacian, symmetric
$\alpha$-stable process, heat kernel, transition density, Green
function, exit time, L\'evy system, intrinsic ultracontractivity,
boundary Harnack principle

\bigskip
\section{Introduction}

It is well-known that, for a second order elliptic differential
operator $\sL$ on $\bR^d$ satisfying some natural conditions, there
is a diffusion process $X$ on $\bR^d$ with $\sL$ as its
infinitesimal generator. The fundamental solution $p(t, x, y)$ of
$\partial_tu=\sL u$ (also called the heat kernel of $\sL$) is the
transition density function of $X$. Thus obtaining sharp two-sided
estimates for $p(t, x, y)$ is a fundamental problem in both analysis
and probability theory.  Such relationship is also true for a large
class of Markov processes with discontinuous sample paths, which
constitute an important family of stochastic processes in
probability theory. They have been widely used in various
applications.

One of the most important and most widely used family of Markov
processes is the family of (rotationally) symmetric $\alpha$-stable
processes on $\bR^d$, $0<\alpha\le 2$. A symmetric $\alpha$-stable
process $X=\{X_t, t\geq 0, \P_x, x\in \bR^d\}$ on $\bR^d$ is a
L\'evy process such that
$$
\E_x \left[ e^{i\xi\cdot(X_t-X_0)} \right]\,=\,e^{-t|\xi|^{\alpha}}
\qquad \hbox{for every } x\in \bR^d \hbox{ and }  \xi\in \bR^d.
$$
When $\alpha=2$, $X$ is a Brownian motion on $\bR^d$ whose
infinitesimal generator is the Laplacian $\Delta$. When $0<\alpha <
2$, the infinitesimal generator of a symmetric $\alpha$-stable
process $X$ in $\bR^d$ is the fractional Laplacian $\Delta^{\alpha
/2}$, which is a prototype of nonlocal operators. The fractional
Laplacian can be written in the form
\begin{equation}\label{e:1.1}
\Delta^{\alpha /2} u(x)
= \,\sA (d, -\alpha) \, \lim_{\eps
\downarrow 0}\int_{\{y\in \bR^d: \, |y-x|>\eps\}} (u(y)-u(x))
\frac{dy}{|x-y|^{d+\alpha}}
\end{equation}
for some constant $ {\cal A}(d, - \alpha):=
\alpha2^{\alpha-1}\pi^{-d/2} \Gamma(\frac{d+\alpha}2)
\Gamma(1-\frac{\alpha}2)^{-1}. $
Here and in the sequel, we use $:=$ as a way of definition.
Here $\Gamma$ is the Gamma function
defined by $\Gamma(\lambda):= \int^{\infty}_0 t^{\lambda-1}
e^{-t}dt$ for every $\lambda > 0$.

Two-sided heat kernel estimates for diffusions in $\bR^d$ have a
long history and many beautiful results have been established. See
\cite{D1, Deb} and the references therein. But, due to the
complication near the boundary, two-sided estimates for the
transition density functions of killed diffusions in a domain $D$
(equivalently, the Dirichlet heat kernels) have been established
only recently. See \cite{D2, Deb, DS} for upper bound estimates and
\cite{Zq3} for the lower bound estimates of the Dirichlet heat
kernels in bounded $C^{1,1}$ domains. In a recent paper \cite{CKS},
we succeeded in establishing sharp two-sided estimates for the heat
kernel of the fractional Laplacian $\Delta^{\alpha /2}$ with zero
exterior condition on $D^c$ (or equivalently, the transition density
function of the killed $\alpha$-stable process) in any $C^{1, 1}$
open set.

The approach developed in \cite{CKS} can be adapted to establish heat
kernel estimates  of other jump processes in open subsets of $\R^d$.
In \cite{CKS1}, the ideas of \cite{CKS} were adapted to establish
two-sided heat kernel estimates of censored stable processes in
$C^{1, 1}$ open subsets of $\R^d$. One of the main tools used in
\cite{CKS1} is the boundary Harnack principle established in
\cite{BBC} and \cite{G}.

In \cite{CKS2} the ideas of \cite{CKS} were adapted to establish
two-sided heat kernel estimates of relativistic stable processes in
$C^{1, 1}$ open subsets of $\R^d$. One of main facts we used in
\cite{CKS2} is that relativistic stable processes can be regarded as
perturbations of symmetric stable processes  in bounded open sets
and therefore the Green functions of killed relativistic stable
processes in bounded open sets are comparable to the Green functions
of killed stable processes in the same open sets.

The goal of this paper is to prove sharp two-sided estimates for the
independent sum of an $\alpha$-stable process and a $\beta$-stable
process, $0<\beta<\alpha<2$, in $C^{1, 1}$ open subsets of $\R^d$.
 Note that these 
processes can  not be  obtained   from symmetric stable
processes through a combination of Girsanov transform and
Feynman-Kac transform.
 So 
the method of \cite{CKS2} can not be
used to establish the comparability of the Green functions of these
processes and the Green functions of symmetric stable processes in
bounded open sets. Since the differences of the L\'evy measures of
these processes and those of symmetric stable processes have
infinite total mass, the method of \cite{R} and \cite{GR} also could
not be used to establish the comparability of the Green functions of
these processes and the Green functions of symmetric stable
processes in bounded open sets.
The approach of this paper will be described in the second paragraph below
after the statement of  Corollary \ref{C:1.2}.

 Let us first 
recall some basic facts about the independent sum of stable
processes and state our main result.

Throughout the remainder of this paper, we assume that $d \ge 1 $
and $0<\beta<\alpha<2$. The Euclidean distance between $x$ and $y$
will be denoted as $|x-y|$. We will use $B(x, r)$ to denote the open
ball centered at $x\in \bR^d$ with radius $r>0$

Suppose $X$ is a symmetric $\alpha$-stable process and $Y$ is a
symmetric $\beta$-stable process on $\R^d$ and that $X$ and $Y$ are
independent. For any $a \ge 0$, we define $X^a$ by $X_t^a:=X_t+ a
Y_t$. We will call the process $X^a$ the independent sum of the
symmetric $\alpha$-stable process  $X$ and the symmetric
$\beta$-stable process  $Y$ with weight $a$. The infinitesimal
generator of $X^a$ is $\Delta^{\alpha/2}+a^\beta \Delta^{\beta/2}$.
Let $p^a(t, x, y)$ denote the transition density of $X^a$ (or
equivalently the heat kernel of $\Delta^{\alpha/2}+ a^\beta
\Delta^{\beta/2}$) with respect to the Lebesgue measure on $\R^d$.
We will use $p(t, x, y)=p^0(t, x, y)$ to denote the transition
density of $X=X^0$ . Recently it is proven in \cite{CK2} that
\begin{eqnarray}
p^1 (t, x, y)\asymp \left(  t^{-d/\alpha} \wedge t^{-d/\beta}
\right)\wedge \left( \frac{t}{|x-y|^{d+\alpha}} +
\frac{t}{|x-y|^{d+\beta}} \right) \quad \mbox{on } (0, \infty)\times
\bR^d \times \bR^d. \label{e:1.0}
\end{eqnarray}
Here and in the sequel,
for $a, b\in \bR$, $a\wedge b:=\min \{a,
b\}$ and $a\vee b:=\max\{a, b\}$;
for any two positive functions $f$ and $g$,
$f\asymp g$ means that there is a positive constant $c\geq 1$
so that $c^{-1}\, g \leq f \leq c\, g$ on their common domain of
definition.

For every open subset $D\subset \bR^d$, we denote by  $X^{a,D}$  the
subprocess of $X^a$ killed upon leaving $D$. The infinitesimal
generator of $X^{a,D}$ is $(\Delta^{\alpha/2} + a^\beta
\Delta^{\beta/2})|_D$, the sum of two fractional Laplacians in $D$
with zero exterior condition. It is known (see \cite{CK2}) that
$X^{a,D}$ has a H\"older continuous transition density $p^a_D(t, x,
y)$ with respect to the Lebesgue measure.

Unlike the case of  the symmetric $\alpha$-stable process $X:=X^0$,
$X^a$ does not have the stable  scaling for $a>0$. Instead, the
following approximate scaling property is true and will be used
several times in the rest of this paper: If $\{X^{a, D}_t, t\geq
0\}$ is the subprocess  of $X^a$  killed upon leaving $D$, then $
\{\lambda^{-1} X^{a, D}_{\lambda^\alpha t}, t\geq 0\}$ is the
subprocess of  $\{X^{a \lambda^{(\alpha-\beta)/\beta}}_t, t\geq 0\}$
killed upon leaving $\lambda^{-1}D$. So for any $\lambda>0$, we have
\begin{equation}\label{e:scaling}
p^{a\lambda^{(\alpha-\beta)/\beta}}_{\lambda^{-1}D} ( t,  x, y) =
\lambda^d p^{a}_D (\lambda^{\alpha}t, \lambda x, \lambda y) \qquad
\hbox{for } t>0 \hbox{ and } x, y \in \lambda^{-1} D.
\end{equation}
In particular, letting $a=1$, $\lambda= a^{\beta/(\alpha -\beta)}$
and $D=\R^d$, we get
$$
p^a(t,x,y)= a^{\frac{\beta d}{\alpha-\beta}}p^1(a^{\frac{\alpha
\beta}{\alpha-\beta}} t ,  a^{\frac{ \beta}{\alpha-\beta}} x ,
a^{\frac{ \beta}{\alpha-\beta}} y) \qquad \hbox{for } t>0
\hbox{ and } x, y\in \R^d.
$$
So we deduce from \eqref{e:1.0} that for any $M>0$ there exists
a constants $C>1$ depending only on $d, \alpha, \beta$
 and $M$ such that for any $a\in (0, M]$ and $(t,x,y) \in (0,
\infty)\times \bR^d \times \bR^d$
\begin{eqnarray}
C^{-1} f^a(t, x, y) \leq p^a (t, x, y)\leq C f^a(t, x, y),
\label{e:1.4}
\end{eqnarray}
where
$$
f^a(t, x, y):=  \left(   (a^\beta t )^{-d/\beta}\wedge
t^{-d/\alpha}\right)\wedge \left(\frac{t}{|x-y|^{d+\alpha}} +
\frac{a^{\beta}t}{|x-y|^{d+\beta}}\right).
$$

\medskip

The purpose of this paper is to establish the following two-sided
sharp estimates on $p^a_D(t, x, y)$ in Theorem \ref{t:main} for
every $t>0$. To state this theorem,  we first recall that an open
set $D$ in $\bR^d$ (when $d\ge 2$) is said to be a (uniform)
$C^{1,1}$ open set if there exist a localization radius $R_0>0$ and
a constant $\Lambda_0>0$ such that for every $z\in\partial D$, there
exist a $C^{1,1}$-function $\phi=\phi_z: \bR^{d-1}\to \bR$
satisfying $\phi (0)= \nabla\phi (0)=0$, $\| \nabla \phi  \|_\infty
\leq \Lambda_0$, $| \nabla \phi (x)-\nabla \phi (z)| \leq \Lambda_0
|x-z|$, and an orthonormal coordinate system $CS_z$ with its origin
at $z$ such that
$$
B(z, R_0)\cap D=\{ y= (\wt y, \, y_d) \mbox{ in } CS_z: |y|< R_0,
y_d
> \phi (\wt y) \}.
$$
The pair $(R_0, \Lambda_0)$ is called the characteristics of the
$C^{1,1}$ open set $D$. Note that a $C^{1,1}$ open set $D$ with
characteristics $(R_0, \Lambda_0)$ can be unbounded and
disconnected; the distance between two distinct components of $D$ is
at least $R_0$. It is well known that any $C^{1, 1}$ open set $D$
satisfies both the {\it uniform interior ball condition} and the
{\it uniform exterior ball condition}: there exists $r_0< R_0$ such
that for every $x\in D$ with $\delta_{\partial D}(x)< r_0$ and $y\in
\bR^d \setminus \overline D$ with $\delta_{\partial D}(y)<r_0$,
there are $z_x, z_y\in \partial D$ so that $|x-z_x|=\delta_{\partial
D}(x)$, $|y-z_y|=\delta_{\partial D}(y)$ and that $B(x_0,
r_0)\subset D$ and $B(y_0, r_0)\subset \bR^d \setminus \overline D$
for $x_0=z_x+r_0(x-z_x)/|x-z_x|$ and $y_0=z_y+r_0(y-z_y)/|y-z_y|$.
By a $C^{1,1}$ open set in $\bR$ we mean an open set which can be
written as the union of disjoint intervals so that the minimum of
the lengths of all these intervals is positive and the minimum of
the distances between these intervals is a positive constant $R_0$.

\begin{thm}\label{t:main}
Suppose $M>0$. Let $D$ be a $C^{1,1}$ open
 subset of $\bR^d$ and $\delta_D(x)$ the Euclidean distance between
$x$ and $D^c$.
\begin{description}
\item{\rm (i)}
For every $T>0$, there is a positive constant
$C_1=C_1(D, M, \alpha, \beta, T)\geq 1$
 such that for every $a \in (0, M]$,
$$
C_1^{-1}\,  f_D^a(t, x, y) \leq p^a_D(t, x, y) \leq C_1\, f_D^a(t,
x, y),
$$
where
$$ f_D^a(t, x, y):= \left( 1\wedge
\frac{\delta_D(x)^{\alpha/2}}{\sqrt{t}}\right) \left( 1\wedge
\frac{\delta_D(y)^{\alpha/2}}{\sqrt{t}}\right) \left( t^{-d/\alpha}
\wedge
\left(\frac{t}{|x-y|^{d+\alpha}}+\frac{a^{\beta}t}{|x-y|^{d+\beta}}
\right)\right).
$$
\item{\rm (ii)} Suppose in addition that $D$ is bounded.
For every $T>0$, there is a constant
$C_2=C_2(D, M, \alpha, \beta, T)\geq 1$
 so that for every $a \in (0, M]$ and $(t, x, y)\in [T, \infty)\times
D\times D$,
$$
C_2^{-1}\, e^{-\lambda_1 t}\, \delta_D (x)^{\alpha/2}\, \delta_D (y)^{\alpha/2}
\,\leq\,  p^a_D(t, x, y) \,\leq\, C_2\, e^{-\lambda_1 t}\, \delta_D (x)^{\alpha/2}
\,\delta_D (y)^{\alpha/2} ,
$$
where $\lambda_1>0$ is the smallest eigenvalue of
$-(\Delta^{\alpha/2} + a^\beta \Delta^{\beta/2})|_D$.
\end{description}
\end{thm}

Letting $a\to 0$, Theorem \ref{t:main} recovers the heat kernel
estimates for symmetric $\alpha$-stable processes obtained in
\cite{CKS}. By integrating the two-sided heat kernel estimates in
Theorem \ref{t:main} with respect to $t$, we obtain the following
estimates on the Green function $ G^a_D(x, y):=\int_0^\infty
p^a_D(t, x, y)dt$, which mean that, for bounded $C^{1, 1}$ domains
$D$, $G^a_D$ and $G^0_D$ are comparable, see \cite{CS1} and
\cite{Ku1}.
To the best of our knowledge, the Green function estimates in the
corollary below are new.

\begin{cor}\label{C:1.2}
Suppose $M>0$. For any bounded $C^{1,1}$-open set $D$ in $\bR^d$,
there is a constant  $C_3=C_3 (D, M, \alpha, \beta)\geq 1$  so that
for every $a \in (0, M]$,
$$
C_3^{-1}\, g_D(x, y) \leq G_D^a(x, y) \leq C_3 \, g_D(x, y)
\qquad \hbox{for } x, y \in D,
$$
where
\begin{equation}\label{e:G}
g_D(x, y):= \begin{cases}
\left(1\wedge \frac{  \delta_D(x)^{\alpha/2}
\delta_D(y)^{\alpha/2}}{ |x-y|^{\alpha}}
\right)\, |x-y|^{\alpha -d}
\qquad &\hbox{when } d>\alpha ,  \\
\log \left( 1+ \frac{  \delta_D(x)^{\alpha/2} \delta_D
(y)^{\alpha/2}}{ |x-y|^{\alpha}}\right)   &\hbox{when } d=1=\alpha , \\
\big( \delta_D(x)  \delta_D (y)\big)^{(\alpha-1)/2} \wedge \frac{
\delta_D(x)^{\alpha/2} \delta_D (y)^{\alpha/2}}{ |x-y|}  &\hbox{when
} d=1<\alpha .
\end{cases}
\end{equation}
\end{cor}

Theorem \ref{t:main}(i) will be established through Theorems
\ref{t:ub} and \ref{t:lb}, which give the upper bound and lower
bound estimates, respectively. Theorem \ref{t:main}(ii) is a
consequence of  the intrinsic ultracontractivity of $X^a$ in  a
bounded open set and the continuity of eigenvalues proved in
\cite{CS8}.  In fact,  the upper bound estimates in both Theorem
\ref{t:main} and Corollary \ref{C:1.2} hold for any
 open set $D$ with
a weak version of the {\it uniform exterior ball condition} in place
of the $C^{1,1}$ condition, while the lower bound estimates in  both
Theorem \ref{t:main} and Corollary \ref{C:1.2} hold for any
open set
$D$ with the {\it uniform interior ball condition} in place of the
$C^{1,1}$ condition (see Theorems \ref{t:ub} and \ref{t:lb}, and the
proofs for Theorem \ref{t:main}(ii) and Corollary \ref{C:1.2}).

Although we follow the main ideas we developed in \cite{CKS}, there
are several new difficulties in  obtaining two-sided Dirichlet heat
kernel estimates for $X^a$: Even though the boundary Harnack
principle has been extended in \cite{KSV} to a large class of pure
jump L\'evy processes including $X^a$, the explicit decay rate of
harmonic functions of $X^a$ near the boundary of $D$ was unknown.
Instead, following the approach in \cite{CKSV}, we establish
necessary estimates using   suitably chosen subharmonic and
superharmonic functions of the process $X^a$. As in \cite{CKSV}, we
need to use finite range (or truncated) symmetric $\beta$-stable
process $\wh Y^\lambda$ obtained from $Y$ by suppressing all its
jumps of size larger than $\lambda$. The infinitesimal generator of
$\wh Y^\lambda$ is
\begin{equation}\label{e:1.7}
\wh \Delta^{\beta/2}_{\lambda} u(x) := {\cal A}(d, -\alpha)\,\lim_{\eps
\downarrow 0}\int_{\{y\in \bR^d: \, \eps<|y-x|\leq \lambda \}} (u(y)-u(x))
\frac{dy}{|x-y|^{d+\beta}}.
\end{equation}
When $\lambda=1$, we will simply denote $\wh
\Delta^{\beta/2}_{\lambda}$ by $\wh \Delta^{\beta/2}$. We first
establish the desired estimates for the L\'evy process $\wh
X^a:=X+a\wh Y^{1/a}$. The infinitesimal generator of $\wh X^a$ is
$\Delta^{\alpha/2}+a^\beta \wh \Delta^{\beta /2}$. The desired
estimates for $X^a=X+aY$ can then be obtained by adding back those
jumps of $Y$ of size larger than $1/a$. To obtain the lower bound of
$p^a(t,x,y)$,  we use the Dirichlet heat kernel estimate for the
fractional Laplacian in \cite{CKS} and a comparison of the killed
subordinate stable process with the subordinate killed stable
process where we will use some of the results obtained in
\cite{SV08}.

We like to point out that unlike \cite{CKS} the boundary Harnack
principle for $X^a$ is not used in this paper, which indicates that
it might be possible to obtain sharp heat kernel estimate  for
processes for which the boundary Harnack principle fails.

As a consequence of Corollary \ref{C:1.2}, we have the following
uniform boundary Harnack principle with explicit decay rate.

\begin{thm}\label{t:bhp}
Suppose that $M>0$. For any $C^{1, 1}$ open set $D$ in $\bR^d$ with
the characteristics $(R_0, \Lambda)$, there exists a positive
constant $C_4=C_4(\alpha, \beta, d, \Lambda, R_0, M)\geq 1$ such
that for $a \in [0, M]$, $r \in (0, R_0]$, $Q\in \partial D$ and any
nonnegative function $u$ in $\R^d$ that is harmonic in $D \cap B(Q,
r)$ with respect to $X^{a}$ and vanishes continuously on $ D^c \cap
B(Q, r)$, we have
\begin{equation}\label{e:bhp_m}
\frac{u(x)}{u(y)}\,\le
C_4\,\frac{\delta^{\alpha/2}_D(x)}{\delta^{\alpha/2}_D(y)} \qquad
\hbox{for every } x, y\in  D \cap B(Q, r/2).
\end{equation}
\end{thm}

Throughout this paper, we will use capital letters $C_1, C_2, \dots$
to denote constants in the statements of results, and their labeling
will be fixed. The lower case constants $c_1, c_2, \dots$ will
denote generic constants used in proofs, whose exact values are not
important and can change from one appearance to another. The
labeling of the lower case constants starts anew in every proof.
The dependence of the constants  on dimension $d$ may not be
mentioned explicitly. For every function $f$, let $f^+:=f \vee 0$.
We will use $\partial$ to denote a cemetery point and for every
function $f$, we extend its definition to $\partial$ by setting
$f(\partial )=0$. We will use $dx$ to denote the Lebesgue measure in
$\bR^d$.
For a Borel set $A\subset \bR^d$, we also use $|A|$ to denote its
Lebesgue measure.

\section{Upper bound estimate}

Throughout this section we assume that $D$ is an open set satisfying
the uniform exterior ball condition with radius $r_0>0$ in the
following sense: for every $z\in \partial D$ and $r\in (0, r_0)$,
there is a ball $B^z$ of radius $r$ such that $B^z\subset \bR^d
\setminus \overline D$ and $\partial B^z \cap \partial D=\{z\}$. The
goal of this section is to establish the upper bound for the
transition density (heat kernel) $p^a_D(t, x, y)$. One of the main
difficulties of getting the  upper bound for $p^a_D(t, x, y)$ is to
obtaining the correct boundary decay rate.

Recall that $\Delta^{\alpha/2}$ and $\wh \Delta^{\beta/2}_\lambda$
are defined by \eqref{e:1.1} and \eqref{e:1.7}. The next two lemmas
can be proved by direct computation, whose proofs
can be found in \cite{G} and \cite{CKSV}, respectively.

For $p>0$, let  $w_p(x):=(x_1^+)^{p}$.

\begin{lemma}\label{this}
For any $x\in (0, \infty)\times \R^{d-1}$, we have
\begin{align}\label{lag1}
\Delta^{\alpha/2} w_{\alpha/2} (x) \ = \ &0.
\end{align}
Moreover, for every $p \in (\alpha/2, \alpha)$, there is a positive
constant $C_5=C_5(d, \alpha, p)$ such that for every $x\in (0,
\infty)\times \R^{d-1}$
\begin{align}\label{lag}
\Delta^{\alpha/2} w_p(x)=&C_5 \, x_1^{p-\alpha}.
\end{align}
\end{lemma}

\begin{lemma}\label{L:2.8}
There are constants $R_* \in (0, 1)$, $C_6> C_7>0$ depending on $p$,
$d$ and $\alpha$ only such that for every $x \in (0,R_*]\times
\R^{d-1}$
\begin{equation}\label{eqn:2.a}
C_7 x_1^{p-\alpha}  \leq \wh \Delta^{\alpha/2} w_p(x)  \leq C_6
x_1^{p-\alpha}  \qquad \hbox{for } \alpha/2 < p < \alpha,
\end{equation}
\begin{equation}\label{eqn:2.b}
| \wh \Delta^{\alpha/2} w_p(x) |  \leq C_6 \, | \log x_1 |
 \qquad \hbox{for } p = \alpha
\end{equation}
and
\begin{equation}\label{eqn:2.c}
| \wh \Delta^{\alpha/2} w_p(x) |  \leq C_6
 \qquad \hbox{for } p > \alpha.
\end{equation}
\end{lemma}

\medskip

In the remainder of this paper, $R_*$ will always stand for the
constant in Lemma \ref{L:2.8}. The following result and its proof
are similar to Lemma 3.2 of \cite{CKSV} and the proof there. For
reader's convenience, we spell out the details of the proof here.

\begin{lemma}\label{L:Main}
Assume that $r_1 \in (0, 1/2] $ and $p \ge \frac{\alpha}2$. Let
$\delta_1:=R_* \wedge (r_1/4)$, $U:=\left\{z\in \bR^d: \, r_1<|z| <
3r_1/2 \right\}$ and
$$
h_p(y):=\left(y_d-\sqrt{r_1^2-|\tilde{y}|^2}\right)^p{\bf 1}_{U \cap
\{z_d >0, |\widetilde{z}| < r_1/2 \}}(y).
$$
Then there exist $C_i=C_i(\alpha, p, r_1)>0$, $i=8, \cdots, 12$,
such that

\begin{description}
\item {\rm (i)}
when $p\in (\alpha/2, \, \alpha)$, we have for all $y \in
\left\{z\in \bR^d: z_d >0, \, r_1<|z| < r_1+\delta_1,
|\widetilde{z}| <r_1/4\right\}$,
 \bee\label{e:h2}
C_8 \left(y_d-\sqrt{r_1^2-|\widetilde{y}|^2}\right)^{p-\alpha}  \leq
\wh \Delta^{\alpha/2} h_p(y)  \le C_9
\left(y_d-\sqrt{r_1^2-|\widetilde{y}|^2}\right)^{p-\alpha}
 \eee
and
 \bee\label{e:h1}
C_8\left(y_d-\sqrt{r_1^2-|\widetilde{y}|^2}\right)^{p-\alpha} \le
\Delta^{\alpha/2}h_p(y)  \le C_9
\left(y_d-\sqrt{r_1^2-|\widetilde{y}|^2}\right)^{p-\alpha};
 \eee

\item{\rm (ii)}
when $p>\alpha$, we have
 \bee\label{e:h0}
|\wh \Delta^{\alpha/2} h_p(y)|  \le C_{10} \quad \text{ for all } y
\in \left\{z\in \bR^d:z_d >0,  \, r_1<|z| < r_1+\delta_1,
|\widetilde{z}| <r_1/4\right\};
 \eee

\item {\rm (iii)}
when $p=\alpha/2$, we have
 \bee\label{e:h3}
|\Delta^{\alpha/2}h_{\alpha/2}(y)|  \le C_{11}   \quad \text{ for
all } y \in \left\{z\in \bR^d: z_d >0, \, r_1<|z| < r_1+\delta,
|\widetilde{z}| <r_1/4\right\};
 \eee

\item {\rm (iv)}
when $p=\alpha$, we have for every  $y \in \left\{z\in \bR^d: z_d
>0, \, r_1<|z| < r_1+\delta, |\widetilde{z}| <r_1/4\right\}$,
 \bee\label{e:h4}
|\wh \Delta^{\alpha/2}h_{\alpha/2}(y)|  \le C_{12} \left|\log
\left(y_d-\sqrt{r_1^2-|\widetilde{y}|^2}\right)\right|  .
 \eee
\end{description}
\end{lemma}

\pf
Let
$$
\Gamma (\tilde{y}):= \sqrt{r_1^2-|\tilde{y}|^2} \quad \text{and}
\quad \underline{h}(y):= y_d- \Gamma (\tilde{y}), \quad y \in U.
$$

Fix $x\in \left\{z\in \bR^d:z_d >0,  \, r_1<|z| < r_1 + R_* \wedge
(r_1/8) , |\widetilde{z}| <r_1/4\right\}$ and choose a point
$x_0\in\partial B_+(0, r_1):=\{z_d >0, |z|=r_1 \}$ satisfying
$\widetilde{x}=\widetilde{x_0}$. Denote by $\overrightarrow{n}(x_0)$
the inward unit normal  vector at $x_0$ for the exterior ball $B(0,
r_1)^c$ and set $\Phi(y)=\langle
y-x_0,\overrightarrow{n}(x_0)\rangle$ for $y\in \bR^d$.
$\Pi=\{y:\Phi(y)=0\}$ is the plane tangent to $\partial B_+(0, r_1)$
at the point $x_0$. Let $\Gamma^*:\widetilde{x}
\in\bR^{d-1}\rightarrow \bR$ be the function describing the
hyperplane $\Pi$, that is, $
\langle(\widetilde{x},\Gamma^*(\widetilde{x}))-x_0,\overrightarrow{n}
(x_0)\rangle=0. $ We also let
\begin{align*}
E&\,:=\,\{y=(\widetilde{y},y_d): y\in U,\ \
|y-x|<r_1/4 \} \,, \\
A&\,:=\,\{y: \Gamma^*(\widetilde{y})>y_d>\Gamma(\widetilde{y}),\
|y-x|<r_1/4\}
\end{align*}
and $\overline{h}(y):=(y_d-\Gamma^*(\widetilde{y})){\bf 1}_{\{ y_d >
\Gamma^*(\widetilde{y})\}}(y)$ for $y\in \bR^d$. Since
$\nabla\Gamma(\widetilde{x})-\nabla\Gamma^*(\widetilde{x})=0$, by
the mean value theorem
\begin{align}\label{hhhh}
|\overline{h}({y})-\underline{h}({y})|\leq|\Gamma(\widetilde{y}) -
\Gamma^*(\widetilde{y})|\leq
\Lambda|\widetilde{y}-\widetilde{x}|^{2},\ \ \ \ \ y\in E.
\end{align}
Let $\delta_{_{\Pi}}(y)={ \rm dist}(y,\Pi)$  for $y\in \bR^d$ and
${U}_{\Gamma^*}=\{y\in \bR^d:y_d>\Gamma^*(\widetilde{y})\}$.
Let $ b_x:=\sqrt{1+|\nabla\Gamma({\widetilde{x}})|^2}$ and
$$
h_{x,p}(y):=(\overline{h}(y))^p.
$$
Note that $h_{x,p}(x)=h_p(x)$ and $B(x,r_1/4)\cap U \subset E$.
Since $\overline{h}(y)=b_x\delta_{_{\Pi}}(y)$ on ${D}_{\Gamma^*}$,
by Lemma \ref{this},
\begin{align}\label{forgqq0}
\Delta^{\alpha/2} h_{x,\alpha/2}(x)\,=\, 0
\end{align}
and, if $\alpha/2 < p <\alpha$,
\begin{align}\label{forgqq}
\Delta^{\alpha/2} h_{x,p}(x)\,= \,c_1 \,b_x^{p}\delta_{_{\Pi}}
^{p-\alpha}(x)  \,= \,c_1 \,b_x^{\alpha}  (\underline{h}(x))^{p-\alpha}
\end{align}
for some $c_1>0$.
By Lemma \ref{L:2.8}, there are constants $c_i>0$, $i=2
\dots 6$, such that for $y\in {D}_{\Gamma^*}$ and
$\delta_{_{\Pi}}(y)< R_*$, when $\alpha/2 < p <\alpha$,
\begin{align}
c_2 (\underline{h}(x))^{p-\alpha} \leq  c_3 \,b_x^p
\delta_{_{\Pi}}^{p-\alpha}(x) \leq \wh \Delta^{\alpha/2} h_{x,p}(x)
\,=\, b_x^p \wh \Delta^{\alpha/2} (\delta_{_{\Pi}}(x))^p \,\le\,c_4
\,b_x^p \delta_{_{\Pi}}^{p-\alpha}(x)  \le c_5
(\underline{h}(x))^{p-\alpha},  \label{e:ee1}
\end{align}
when $p> \alpha$,
\begin{align}
|\wh \Delta^{\alpha/2} h_{x,p}(x)|\,=\,b_x^p | \wh
\Delta^{\alpha/2} (\delta_{_{\Pi}}(x))^p |\,\le\,c_6
 \label{e:ee2}
\end{align}
and when $p=\alpha$,
\begin{equation}\label{e:ee3}
|\wh \Delta^{\alpha/2} h_{x,p}(x)|\leq c_6 \, | \log
(\underline{h}(x)) |.
\end{equation}
Note that
\begin{align}
| \wh \Delta^{\alpha/2} (h_p-h_{x,p})(x)| & = \sA (d, -\alpha)
\bigg|\lim_{\varepsilon\downarrow 0}\int_{
\{1\geq|y-x|>\varepsilon\}}\frac{(h_p(y)-h_{p,x}(y))}{|x-y|^{d+\alpha}}\
dy\bigg|\nonumber\\
&\leq \, \sA (d, -\alpha)  \bigg| \int_{ \{1\geq
|y-x|>r_1/4\}}\frac{(h_p(y)-h_{p,x}(y))}{|x-y|^{d+\alpha}}\
dy\bigg|\\
&\qquad +\sA (d, -\alpha)   \lim_{\varepsilon\downarrow 0}\int_{
\{r_1/4\geq |y-x|>\varepsilon\}}\frac{|h_p(y)-h_{p,x}(y)|}{|x-y|
^{d+\alpha}}\ dy \nonumber\\
&\,\leq \ c_7 +\sA (d, -\alpha)\int_{A} \frac{h_p(y)+h_{p,x}(y)}
{|x-y|^{d+\alpha}}\ dy +\sA (d, -\alpha) \int_{E}
\frac{|h_p(y)-h_{p,x}(y)|}{|x-y|^{d+\alpha}}  \nonumber\\
&=:\,c_7+I_1+I_2
 \label{e:I}
\end{align}
and, similarly,
\begin{align}
&| \Delta^{\alpha/2} (h_p-h_{x,p})(x)| \nonumber  \\
\leq & \sA (d,-\alpha)
\bigg| \int_{
\{|y-x|>r_1/4\}}\frac{(h_p(y)-h_{p,x}(y))}{|x-y|^{d+\alpha}}\
dy\bigg|\nonumber\\
&+\sA (d, -\alpha)\int_{A} \frac{h_p(y)+h_{p,x}(y)}
{|x-y|^{d+\alpha}}\ dy +\sA (d, -\alpha) \int_{E}
\frac{|h_p(y)-h_{p,x}(y)|}{|x-y|^{d+\alpha}}
=:I_3+I_1+I_2. \label{e:I2}
\end{align}
Since for $y\in B(x,r_1/4)^c$,
$$
|h_{x,p}(y)-h_{x,p}(x)| \le c_8 |y-x|^{p} \quad \text{ and } \quad
|h_p(y)| \le c_8
$$
and  $h_p(y)=0$ for $|\widetilde{y}|>r_1/2$, for $\alpha/2 \le p <
\alpha$ we get
\begin{align}
I_3 \le & \sA (d, -\alpha)
\int_{B(x,r_1/4)^c}\frac{|h_{x,p}(y)-h_{x,p}(x)|}{|x-y|^{d+\alpha}}\
dy+ \sA (d, -\alpha) \int_{B(x,r_1/4)^c \cap \{|\widetilde{y}|\le
r_1/2\}}\frac{|h_p(y)-h_p(x)|}{|x-y|^{d+\alpha}}
dy\nonumber\\
&+ \sA (d, -\alpha) \bigg|\int_{B(x,r_1/4)^c \cap
\{|\widetilde{y}|>r_1/2\}}\frac{h_p(x)}{|x-y|^{d+\alpha}}\ dy\bigg|
\nonumber\\
\leq&c_9\int_{B(x,r_1/4)^c}\frac{1}{|x-y|^{d+\alpha-p}}
dy+c_9\int_{B(x,r_1/4)^c}\frac{1}{|x-y|^{d+\alpha}} dy \le c_{10}
<\infty. \label{*}
\end{align}
We claim that, if $p \ge \alpha/2$,
 \bee\label{e:I1234}
 I_1\,+\,I_2 \,\leq\, c_{11} \,<\, \infty.
 \eee
Note that for $y \in A$
\begin{align}
|h_{x,p}(y)|+|h_{p}(y)| &\le |y_d-\Gamma^*(\widetilde{y})|^p
+|y_d-\Gamma(\widetilde{y})|^p\le  2|\Gamma(\widetilde{y})-
\Gamma^*(\widetilde{y})|^p \nonumber \\
&\le2 |\Gamma(\widetilde{y})-\Gamma(\widetilde{x})-
\nabla\Gamma(\widetilde{x})\cdot (\widetilde{y}-\widetilde{x})|^p
\le  2 c_{12}^{p}|\widetilde{y}-\widetilde{x}|^{2p} \label{e:dfe2}.
\end{align}
Furthermore, since $|\Gamma(\widetilde{y})-\Gamma^*(\widetilde{y})| \le
c_{13}|\widetilde{y}-\widetilde{x}|^{2}\leq c_{12} r^2$ on
$|y-x|=r$, this together with \eqref{e:dfe2} yields that
\begin{eqnarray*}
I_1 &\leq& c_{14}  \int_{ 0}^{r_1/4}  r^{2p-\alpha-d} \int_{|y-x|=r}{\bf
1}_A(y)\ m_{d-1}(d y)dr\nonumber\\
&= & c_{14} \int_{ 0}^{ r_1/4} r^{2p-\alpha-d} m_{d-1} \big(\{ y:\,
| y-x|=r, \ \Gamma^*(\widetilde{y})>y_d> \Gamma(\widetilde{y}) \} \big) dr\\
&\leq& c_{15} \int_{ 0}^{ r_1/4} r^{2p-\alpha} dr  \,<\, \infty .
\end{eqnarray*}

Note that for $y\in E$
\begin{align}
|h_p(y)-h_{x,p}(y)| \le& c_{16} | (\overline{h}(y))^{p}
-(\underline{h}(y))^{p}| \le c_{17}  (\overline{h}(y))^{(p-1)_-}|
\overline{h}(y) -\underline{h}(y)|, \label{e:KG}
\end{align}
where $(p-1)_-:=(p-1) \wedge 0$. In the last inequality above, we
have used the inequalities
\begin{align*}
|b^{p}-a^{p}|\leq b^{p-1}|b-a| \qquad \hbox{for } a , b>0, \ 0 <p
\le 1
\end{align*}
and
\begin{align*}
|b^{p}-a^{p}|\leq (p+1) |b-a| \qquad \hbox{for } a,b \in (0, 1),  \
p >1 .
\end{align*}
For $y=(\wt y, y_d) \in \R^d$, we use an affine coordinate system $z
= (\wt z, z_d)$ to represent it so that $z_d= y_d-\Gamma^*(\wt y)$
and $\wt z$ is the coordinates in an orthogonal coordinate system
centered at $x_0$ for the $(d-1)$-dimensional hyperplane $\Pi$ for
the point $(\wt y, \Gamma^*(\wt y))$. Denote such an affine
transformation $y\mapsto z$ by $z=\Psi (y)$. It is clear that there
is a constant $c_{18}>1$ so that for every $y\in \R^d$,
$$
c_{18} ^{-1}|\wt y -\wt x| \leq |\wt z|\leq c_{18} |\wt y -\wt x|, \qquad
c_{18}^{-1} |y-x|\leq |\Psi (y)-\Psi (x)| \leq c_{18} |y-x|
$$
and that
$$
\Psi(E)\subset \{z=(\wt z, z_d)\in \R^d: \ |\widetilde{z}|<c_{18} r_1
\hbox{ and }  0< z_d\leq c_{18} r_1\}.
$$
Denote $x_d-\Gamma^*(\wt x)$ by $w$; that is,
 $\Psi (x)=(\wt 0, w)$.
Hence by \eqref{hhhh} and \eqref{e:KG} and applying the transform
$\Psi$,  we have by using polar coordinates for $\wt z$ on the
hyperplane $\Pi$,
\begin{eqnarray*}
I_2 &\leq&c_{19}\int_{ E}\frac{ \overline{h}(y)^{(p-1)_-}\,
|\widetilde{y}-\widetilde{x}|^{2}}{|y-x|^{d+\alpha}}\, dy \leq
c_{19} \int_{\Psi(E)}\ \frac{z_d^{(p-1)_-} \,
|\wt z|^{2}}{|z- (\wt 0, w)|^{d+\alpha}}\ dz\nonumber\\
&\leq & c_{20}\int_{0}^{c_{18}r_1}  z_d^{(p-1)_-}\left(
\int_0^{c_{18}r_1} \frac{r^{d-2}}{(r+|z_d-w|)^{d+\alpha-2}} dr
\right) dz_d\\
&\leq & c_{20}\int_{0}^{c_{18}r_1}  z_d^{(p-1)_-}\left(
\int_0^{c_{18}r_1}
\frac{1}{(r+|z_d-w|)^{\alpha}} dr \right) dz_d\\
&\leq& c_{21}\int_{0}^{c_{18}r_1} z_d^{(p-1)_-}\left(
\frac1{|z_d-w|^{\alpha-1}} -
\frac{1}{(c_{18}r_1+|z_d-w|)^{\alpha-1}} \right)dz_d\\
&<&  c_{22}\int_{0}^{c_{18}r_1}
\frac1{z_d^{(1-p)^+}\,
|z_d-w|^{\alpha-1}} dz_d  \leq  c_{23}<\infty,
\end{eqnarray*}
where all the constants depend on $\alpha$, $p$ and $r_1$. The last
inequality is due to the fact that since $p>0$, $0<\alpha<2$ and
$(1-p)^+ + \alpha -1 =\max\{\alpha-p, \alpha -1\}<1$, by the
dominated convergence theorem, $\phi (w):=
\int_0^{c_{18}r_1}\frac1{z_d^{(1-p)^+}\,  |z_d-w|^{\alpha-1}} dz_d$
is a strictly positive continuous function in $x_d\in [0,
c_{18}r_1]$ and hence is bounded.

Thus we have proved the claim \eqref{e:I1234}. The desired estimates
\eqref{e:h2}-\eqref{e:h4} now follow from
\eqref{forgqq0}-\eqref{e:I1234}. \qed

It is well-known that  $X^1$ has L\'evy intensity
$$
J^1(x, y)=j^1(|x-y|)=\frac{{\cal A}(d, -\alpha)}{|x-y|^{d+\alpha}} +
\frac{{\cal A}(d, -\beta)}{|x-y|^{d+\beta}}.
$$
A scaling argument yields that
$$
J^a(x, y)=j^a(|x-y|)=\frac{{\cal A}(d, -\alpha)}{|x-y|^{d+\alpha}}+
a^{\beta}\frac{{\cal A}(d, -\beta)}{|x-y|^{d+\beta}}.
$$
Put
\begin{align}
\psi^a(r):=1+a^{\beta}\frac{{\cal A}(d, -\beta)}{{\cal A}(d,
-\alpha)}r^{\alpha-\beta}, \quad r\ge 0. \label{e:psi}
\end{align}
Clearly for $a\in (0, M]$ and $r>0$, $ 1 \le \psi^a(r) \le 1+c
M^{\beta} r^{\alpha-\beta}$ and
$$
J^a(x, y)=j^a(|x-y|)=\frac{{\cal A}(d,
-\alpha)}{|x-y|^{d+\alpha}}\psi^a(|x-y|).
$$

The L\'evy intensity gives rise to a L\'evy system for $X^a$, which
describes the jumps of the process $X^a$: for any non-negative
measurable function $f$ on $\bR_+ \times \bR^d\times \bR^d$, $x\in
\bR^d$ and stopping time $T$ (with respect to the filtration of
$X^a$),
\begin{equation}\label{e:levy}
\E_x \Big[\sum_{s\le T} f(s,X^a_{s-}, X_s) \Big]=
\E_x \Big[ \int_0^T \Big( \int_{\bR^d} f(s,X^a_s, y) J^a(X^a_s,y) dy
\Big) ds \Big].
 \end{equation}
(See, for example, \cite[Proof of Lemma 4.7]{CK} and \cite[Appendix A]{CK2}).

For any open set $D\subset \bR^d$, let $\tau^a_D=\tau^a(D): =
\inf\{t>0: \,X^a_t\notin D\}$ denote the first exit time from $D$ by
$X^a$.

The next lemma follows immediately from a special case of
\cite[Proposition 2.10 and Lemma 3.6]{KSV} and the scaling property
of $X^a$.

\begin{lemma}\label{L:2.00}
For any $b, M \in (0,\infty)$, there exists
$C_{13}=C_{13}(M, b, \alpha, \beta)>0$
 such that for every $x_0 \in \bR^d$, $a \in [0,M]$ and $r \in (0,
b]$,
 \bee \label{e:ext}
\E_x\left[\tau^a_{B(x_0,r)}\right]\, \le\, C_{13}\, r^{\alpha/2}
(r-|x-x_0|)^{\alpha/2}, \qquad \hbox{for } x \in B(x_0,r).
 \eee
\end{lemma}

\medskip

For $\lambda >0$, $\wh Y^\lambda=(\wh Y_t^\lambda, \P_x)$ is a
L\'evy process in $\R^d$ such that
$$
\E_x\left[e^{i\xi\cdot(\wh Y^\lambda_t-\wh Y^\lambda_0)}\right] =
e^{-t\psi(\xi)} \quad \quad \mbox{ for every } x\in \R^d
\mbox{ and } \xi\in \R^d,
$$
with
$$
\psi(\xi)=\sA(d, -\beta) \int_{\{|y|\leq
\lambda\}}\frac{1-\cos(\xi\cdot y)}{|y|^{d+\beta}}dy.
$$
In other words, $\wh Y^\lambda$ is a pure jump symmetric  L\'evy
process in $\R^d$ with a  L\'evy density given by ${\cal A}(d,
-\beta) |x|^{-d-\beta}\, 1_{\{|x|\leq \lambda\}}$. For $a>0$,
suppose $\wh Y^{1/a}$
is independent of the symmetric $\alpha$-stable process $X$ on
$\R^d$. Define
$$
\wh X_t^{a}:=X_t+ a\wh Y_t^{1/a}, \qquad t\geq 0 .
$$
We will call the process $\wh X^{a}$ the independent sum of the
symmetric $\alpha$-stable process $X$ and the truncated symmetric
$\beta$-stable process $\wh Y^{1/a}$ with weight $a>0$. The
infinitesimal generator of $\wh X^{a}$ is $\Delta^{\alpha/2} +
a^\beta \wh \Delta^{\beta/2}$.

For any open set $U\subset \bR^d$, let $\wh \tau^{a}_U=\inf\{t>0: \,
\wh X^{a}_t\notin U\}$ be the first exit time from $U$ by $\wh
X^{a}$. The truncated process $\wh X^{a}$  will be used in the proof
of next lemma.

\begin{lemma}\label{L:2.0}
Assume $r_1 \in (0, \frac14]$ and $M>0$.  Let  $U:=\left\{z\in
\bR^d: \, r_1<|z| < 3r_1/2 \right\}$. There are  constants
$C_{14}=C_{14}(r_1, \alpha)>0$ and $C_{15}=C_{15}(r_1, M, \alpha,
\beta)>0$
 such that for every $a \in [0, M]$
 \bee\label{e:L:2.0}
\E_x[\tau^a_U] \,\le\, C_{14} \P_x \left( |X^{a}_{\tau^a_U}| \ge
3r_1/2\right) \,\le\,C_{15}\delta_U(x)^{\alpha/2}, \qquad \hbox{for
} r_1<|x|<5r_1/4.
 \eee
\end{lemma}

\pf The first inequality in \eqref{e:L:2.0} is easy. In fact, by the
L\'evy system \eqref{e:levy} with
$$
f(s, x, y) ={\bf 1}_U(x) {\bf 1}_{\{5r_1<|y|<10r_1\}}(y)
$$
and $T=\tau^a_U$, we have that for $x \in U$
\begin{eqnarray*}
&&\P_x \left( |X^{a}_{\tau^a_U}| \ge 3r_1/2\right) \,\ge\, \P_x
\left( 10r_1> |X^a_{\tau^a_U}|> 5r_1 \right)\\
&& = \E_x \left[ \int_0^{\tau^a_U} \int_{\{ 5r_1<|y|<10r_1\}}
J^a(X^a_s, y)dy ds\right]\\
&& \ge \E_x \left[ \int_0^{\tau^a_U} \int_{\{ 5r_1<|y|<10r_1\}}
\frac{\sA (d, -\alpha)}{|X^a_s-y|^{d+\alpha} } dy ds\right] \ge c_1
\E_x [\tau^a_U] ,
\end{eqnarray*}
where $c_1=c(r_1, \alpha)>0$.

It is enough to prove the second inequality in \eqref{e:L:2.0} for
$r_1<|x|< r_1+\delta$ for some small $\delta>0$. Without loss of
generality, we assume
$\widetilde x =\widetilde 0$
and $x_d >0$.
 Let $p>0$ be
such that $p\not=\beta$ and
$$
\alpha -(\beta/2)<p< (\alpha - (\beta/2)+(\alpha-\beta)/3)\wedge
\alpha.
$$
Note that $ \alpha/2<p<3\alpha/2 - \beta$. Define
$$
h(y):=\left(y_d-\sqrt{r_1^2-|\widetilde{y}|^2}\right)^{\alpha/2}
{\bf 1}_{U\cap \{z_d >0, |\widetilde{z}|<r_1/2\}}(y),
$$
$$
g_p(y):=\left(y_d-\sqrt{r_1^2-|\widetilde{y}|^2}\right)^p{\bf
1}_{U\cap \{z_d >0, |\widetilde{z}|<r_1/2\}}(y),
$$
and let $\phi$ be a smooth function on $\bR^d$ with bounded first
and second partial derivatives such that $ \phi(y)=2^{4+p}
|\widetilde y|^2/r_1^2$ for
$ y \in \{z_d >0,   r_1<|y|<4r_1/5,   |\widetilde z|<r_1/4\}$
and $2^p \le \phi(y)\le  4^p$ if $|\widetilde y| \ge
r_1/2$ or $|y| \ge 3r_1/2. $

Since $r_1 \le 1/4$, it is easy to see that $\|g_p\|_\infty<1$. Now
we define
$$
u(y):=h(y) +  \phi(y) -g_p(y).
$$
By Taylor's expansion with the remainder of order $2$, we get that
for any $a\in (0, M]$ and $y\in \R^d$,
 \bee\label{e:tr2}
\big|(\Delta^{\alpha/2}+a^\beta\wh \Delta^{\beta/2})\phi(y)\big| \, \le\,
\big\|\Delta^{\alpha/2}\phi\big\|_\infty+M^\beta\big\|\wh
\Delta^{\beta/2}\phi\big\|_\infty \, \le \,c_2(\alpha, \beta, M)\,<
\,\infty.
 \eee
Moreover, by \eqref{e:h2}--\eqref{e:h0}, there
 exist $c_3=c_3(\alpha, \beta)>0$ and $\delta_1=\delta_1(\alpha,
\beta)\in (0, r_1/8)$ such that
$$
\Delta^{\alpha/2}g_p(y) \ge  c_3 \delta_D(y)^{p-\alpha} \quad
\text{for } y \in \left\{z\in \bR^d:z_d >0,  \, r_1<|z| <
r_1+\delta_1,  |\widetilde{z}|<r_1/4 \right\}
$$
and
$$
\wh \Delta^{\beta/2}g_p(y) \ge -c_3 \delta_D(y)^{(p-\beta)\wedge 0}
\quad \text{for } y \in \left\{z\in \bR^d: z_d >0, \, r_1<|z| <
r_1+\delta_1 , |\widetilde{z}|<r_1/4\right\}.
$$
Thus there exist $c_4=c_4(\alpha, \beta, M) >0$ and $\delta_2
=\delta_2(\alpha, \beta, M)\in (0, \delta_1)$ such that  for all
$a\in (0, M]$ and $ y \in \left\{z\in \bR^d: z_d >0, \, r_1<|z| <
r_1+\delta_2, |\widetilde{z}|<r_1/4\right\}$,
 \bee\label{e:gp}
(\Delta^{\alpha/2}+ a^\beta \wh \Delta^{\beta/2})g_p(y) \ge c_3
\delta_D(y)^{p-\alpha} -c_3 M^\beta \delta_D(y)^{(p-\beta)\wedge 0}
\ge c_4 \delta_D(y)^{p-\alpha}.
 \eee
Furthermore by \eqref{e:h2} and \eqref{e:h0}--\eqref{e:h4}, there
exist $c_5=c_5(\alpha, \beta, M)>0$ and $\delta_3 =\delta_3(\alpha,
\beta)\in (0, \delta_1)$ such that for all $a\in (0, M]$ and for
every $y \in \left\{z\in \bR^d: z_d >0, \, r_1<|z| < r_1+\delta_3 ,
|\widetilde{z}|<r_1/4\right\},$
\begin{eqnarray}\label{e:gp2}
\left|(\Delta^{\alpha/2}+ a^\beta \wh \Delta^{\beta/2}) h(y) \right|
& \le&  \big| \Delta^{\alpha/2}h(y) \big| + M^\beta \big| \wh
\Delta^{\beta/2} h(y)\big| \nonumber \\
& \le & \begin{cases}
c_5 + c_5\delta_D(y)^{(\alpha/2-\beta) \wedge 0} &\hbox{if } \beta \not= \alpha/2, \\
c_5 + c_5 | \log \delta_D (y)| & \hbox{if } \beta = \alpha/2 .
\end{cases}
\end{eqnarray}
Since $p-\alpha <\alpha/2-\beta$, by  \eqref{e:tr2}--\eqref{e:gp},
there exists $\delta_4=\delta_4(\alpha, \beta, M) \in (0, \delta_2
\wedge \delta_3)$ such that for all $a\in (0, M]$ and $ y \in
V:=\left\{z\in \bR^d: z_d >0, \, r_1<|z| < r_1+\delta_4 ,
|\widetilde{z}|<r_1/4\right\}$
 \bee\label{e:sh}
(\Delta^{\alpha/2}+a^\beta \wh \Delta^{\beta/2})u(y) \le c_2+c_5 +
 c_5 \left(\delta_D(y)^{(\alpha/2-\beta) \wedge 0} + | \log \delta_D(y)|\right)
  - c_4 \delta_D(y)^{p-\alpha} \le 0.
 \eee

Let $\eta$ be a non-negative smooth radial  function with compact
support in $\R^d$ such that $\eta(x)=0$ for
 $|x|>1$ and $\int_{\R^d} \eta (x) dx=1$. For $k\geq 1$, define
$\eta_k(x)= 2^{kd} \eta (2^k x)$. Set $ u^{(k)}(z):= ( \eta_k*u)(z).
$ As $(\Delta^{\alpha/2}+a^\beta \wh \Delta^{\beta/2})u^{(k)} =
\eta_k* (\Delta^{\alpha/2}+a^\beta \wh \Delta^{\beta/2}) u, $ we
have by \eqref{e:sh} that
$$
\big(\Delta^{\alpha/2}+ a^\beta \wh \Delta^{\beta/2}\big)u^{(k)}\le
0
$$
on $V_k:=\left\{z\in \bR^d: z_d >0, \, r_1+ 2^{-k}<|z| <
r_1+\delta_4 -2^{-k}  \text{ and }  |\widetilde{z}|<r_1/4
-2^{-k}\right\}$. Since $u^{(k)}$ is a bounded smooth function on
$\R^d$ with bounded first and second partial derivatives, by Ito's
formula and the L\'evy system \eqref{e:levy},
\begin{equation}\label{e:cl0}
M^k_t:=u^{(k)}(\wh X^{a}_t)-u^{(k)}(\wh X^{a}_0)-\int_0^t
\left(\Delta^{\alpha/2}+
a^\beta
 \wh \Delta^{\beta/2}\right)u^{(k)}(\wh
X^{a}_s)ds
\end{equation}
is a martingale. Thus it follows from \eqref{e:cl0} that $ t\mapsto
u^{(k)} \big( \wh X^{a}_{t\wedge \wh \tau^{a}_{V_k}}\big) $ is a
bounded supermartingale. Since $V_k$ increases to $V$ and $u$ is
bounded and continuous on $\overline V$, we conclude that
\begin{equation}\label{e:claim1}
t\mapsto u \Big( \wh X^{a}_{t\wedge \wh \tau^{a}_{V}}\Big) \  \hbox{
is a bounded supermartingale}.
\end{equation}

We observe that, since $\phi(x)=0$,
 \bee\label{e:phi0}
u(x) \le \delta_U(x)^{\alpha/2}.
 \eee
We also observe that, since $\phi \ge 2 g_p$ outside of $\left\{z\in
U: \, z_d>0,\, |\widetilde{z}|< r_1/2\right\}$ and
$$
u(y)
\ge  \left(y_d-\sqrt{r_1^2-|\widetilde{y}|^2}\right)^{\alpha/2}
-\left(y_d-\sqrt{r_1^2-|\widetilde{y}|^2}\right)^p > c_6
$$
on $\{ z_d>0,\,   r_1+\delta_4 \le|z| < 3r_1/2, \, |\widetilde{z}|<
r_1/2\},$ we have
\begin{align} \label{e:dddd}
u(y) \ge c_7>0 \quad \text{ for } y\in V^c \setminus
\overline{B(0,r_1)},
\end{align}
where $c_7$ depends on $\delta_4$, $\alpha$, $\beta$ and $r_1$.
Therefore, by \eqref{e:claim1}-\eqref{e:dddd} we get
\begin{align}
\delta_U(x)^{\alpha/2} \ge u(x) \ge \E_{x}\left[u\big( \wh X_{\wh
\tau^a_V}\big)\right] \ge c_7 \P_{x}\left(\wh X^{a}_{\wh \tau^a_V} \in
V^c \setminus \overline{B(0,r_1)}\right) \ge c_7 \P_{x}\left(|\wh
X^{a}_{\wh \tau^a_U}| \ge  3r_1/2\right). \label{e:HA1}
\end{align}

Note that there exist $c_8=c_8(\alpha, d, r_1)>0$ and
$c_{9}=c_{9}(\beta, d, r_1)>0$ such that  for $z\in U$,
$$
\int_{\{ |y| \ge 2r_1 \}} \frac{dy}{|z-y|^{d+\alpha}} \le c_8\int_{\{
2r_1 \le |y| <3r_1\}} \frac{dy}{|z-y|^{d+\alpha} }
$$
and
$$
\int_{\{ |y| \ge 2r_1 \}} \frac{dy}{|z-y|^{d+\beta} } \le
c_{9}\int_{\{ 2r_1 \le |y| <3r_1\}} \frac{dy}{|z-y|^{d+\beta} }.
$$
Thus by \eqref{e:levy}, there exists a positive constant
$c_{10}=c_{10}(d, \alpha, \beta, M)$ such that for any $a\in (0,
M]$,
\begin{align}
\P_{x}\left(| X^{a}_{ \tau^a_U}| \ge  2 r_1\right)\,
&=\, \E_x \left[ \int_0^{\tau^a_U} \int_{\{ |y| \ge 2r_1\}}
J^a(X^a_s, y) dy ds\right]\nn\\
&\le\, c_{10} \E_x \left[ \int_0^{\tau^a_U} \int_{\{ 2r_1 \le |y| <
3r_1\}} J^a(X^a_s, y) dy ds\right]\nn\\
&=\,c_{10}\P_{x}\left(3r_1>| X^{a}_{ \tau^a_U}| \ge  2 r_1\right).
\label{e:HA2}
\end{align}
Since $r_1\leq 1/4$ and the processes $X$ and $Y$ do not jump
simultaneously, we have  by \eqref{e:HA1} that there is a positive
constant $c_{11}=c_{11}(d, \alpha, \beta, M, r_1)$ such that for all
$a\in (0, M]$,
\begin{align*}
 \P_{x}\left(| X^{a}_{ \tau^a_U}|  \ge  3r_1/2\right)\, &\le\,
(c_{10}+1) \P_{x}\left(3r_1>| X^{a}_{ \tau^a_U}| \ge
3r_1/2\right)\\
&=\, (c_{10}+1) \P_{x}\left(3r_1>| \wh X^{a}_{\wh \tau^a_U}|
\ge  3r_1/2\right) \\
& \le\, (c_{10}+1) \P_{x}\left(| \wh X^{a}_{\wh \tau^a_U}| \ge
3r_1/2\right)\,\le\, c_{11} \delta_U(x)^{\alpha/2}.
\end{align*}
\qed

\begin{lemma}\label{L:2.1}
Assume $M>0$ and $r_1 \in (0, \frac14]$. Let $E = \{x\in \bR^d: \,
|x| >r_1\}$. Then  for every $T>0$, there is a constant
$C_{16}=C_{16}( r_1, \alpha, \beta, T, M)>0$
 such that for every  $a \in [0, M]$,
$$
p^a_E(t, x, y) \leq C_{16} \, {\delta_E(x)^{\alpha/2}}J^a(x, y)
\qquad \hbox{for } r_1<|x|<5r_1/4, \ |y|\geq 2r_1 \hbox{ and } t\leq
T.
$$
\end{lemma}

\pf Define $U:=\left\{z\in \bR^d: \, r_1<|z| < 3r_1/2 \right\}$.
Since $X^a$ satisfies the hypothesis ${\bf H}$ in \cite{Sz1}, by
\cite[Theorem 1]{Sz1}, $X^a_{\tau^a_U}\notin \partial U$ with
 probability 1. For $r_1<|x|<5r_1/4$, $ |y|\geq 2r_1$ and $t\in (0,
T]$, it follows from the strong Markov property of $X^a$ and
\eqref{e:levy}  that
\begin{align*}
&p^a_E(t, x, y)\,=\, \E_x\left[ \,p^a_E(t-\tau^a_U, X^a_{\tau^a_U}, y):
\, \tau^a_U < t  \,\right]\nonumber\\
=&\int_0^t  \int_U p_U(s, x, z) \left( \int_{\{ w:
|w|>3r_1/2\}} J^a(z, w)\ p^a_E(t-s, w, y) dw\right) dzds
\nonumber \\
=&  \int_0^t  \int_U p_U(s, x, z) \left( \int_{\{ w: \,
(3r_1/4)+(|y|/2)\geq |w|>3r_1/2\}} J^a(z, w) p^a_E(t-s, w, y)
dw\right) dz ds  \nonumber \\
& + \int_0^t \int_U p_U(s, x, z) \left( \int_{\{ w: \, |w|>
(3r_1/4)+(|y|/2)\}} J^a(z, w) p^a_E(t-s, w, y) dw\right) dz
ds  \nonumber \\
=:& I+II.
\end{align*}
Note that for $|w|\leq  (3r_1/4)+(|y|/2)$,
 \bee\label{e:wy}
|w-y| \geq |y|-|w| \geq \frac12 \left(|y|-\frac{3r_1}2\right)\geq
\frac{|y|}8\geq \frac{|x-y|}{16}.
 \eee
Since $p^a_E (t-s, w, y)\leq p^a (t-s, w, y)$, by \eqref{e:1.4} and
\eqref{e:wy}, there exist constants $c_1=c_1(\alpha, \beta, M)>0$
and $c_2=c_2(\alpha, \beta, M)>0$ such that for $a\in (0, M]$
\begin{align*}
I &\leq  \int_0^t  \int_U p^a_U(s, x, z) \left( \int_{\{ w: \,
(3r_1/4)+(|y|/2)\geq |w|>3r_1/2\}} J^a(z, w) \, c_1T J^a(w, y)
dw\right) dz ds \\
&\leq  c_2T J^a(x, y) \int_0^t  \int_U p^a_U(s, x, z) \left(
\int_{\{ w: \, 3|x-y|/4\geq |w|>3r_1/2\}} J^a(z,
w)  dw\right) dz ds \\
&= c_2T J^a(x, y) \, \P_x \left( 3r_1/2< |X^a_{\tau^a_U}| \leq
3|x-y|/4; \, \tau^a_U
\leq t \right) \\
&\leq  c_2T J^a(x, y) \, \P_x \left( |X^a_{\tau^a_U}| > 3r_1/2
\right).
\end{align*}
By Lemma \ref{L:2.0} , we have for $|x|\in (r_1, 5r_1/4)$,
$$
\P_x \left(   |X^a_{\tau^a_U}|> 3r_1/2 \right) \leq  c_3 \,\delta_U
(x)^{\alpha/2}= c_3\delta_E (x)^{\alpha/2}
$$
for some positive constant $c_3=c_3(r_1, \alpha, \beta, M)$. Thus
 \bee\label{e:2.3}
I\leq c_4 \, (T\vee 1)\, { \delta_E (x)^{\alpha/2}} J^a(x, y)
 \eee
for some positive constant $c_4=c_4(r_1, \alpha, \beta, M)$. On the
other hand, for $z\in U$ and $w\in \bR^d$ with
$|w|>(3r_1/4)+(|y|/2)$,
$$
 |z-w| \geq |w|-|z|
 \geq  \frac12\left(|y|-\frac{3r_1}2\right) \geq
\frac{|y|}8\geq \frac{|x-y|}{16}.
$$
Thus by the symmetry of $p^a_E(t-s, w, y)$ in $(w, y)$, we have
\begin{align*}
II
&\leq  c_5J^a(x, y) \int_0^t \int_U p^a_U(s, x, z) \left( \int_{\{ w:
\,  |w|> (3r_1/4)+(|y|/2)\}}  p^a_E(t-s, y, w)
dw\right) dz ds \\
&\leq   c_5 J^a(x, y) \, \int_0^\infty
 \int_U p^a_U(s, x, z) dz  ds \\
&=  c_5 J^a(x, y) \E_x \left[ \tau^a_{U}\right]\
\leq \ c_6 \, \delta_E (x)^{\alpha/2} J^a(x, y)
\end{align*}
for some positive constants $c_k=c_k (r_1, \alpha, \beta, M)$, $k=5,
6$. In the last inequality, we used Lemma \ref{L:2.0} to deduce that
$\E_x \left[ \tau^a_U \right] \leq c \delta_U(x)^{\alpha/2} = c
\delta_E(x)^{\alpha/2}$ for some positive constant $c=c(r_1, \alpha,
\beta, M)$. This together with \eqref{e:2.3} proves the lemma. \qed

\medskip

\begin{thm}\label{ub11}
Assume that $M>0$ and $D$ is an open set that satisfies the uniform
exterior ball condition with radius $r_0>0$. Then for every $T>0$,
there is a constant $C_{17} = C_{17}(r_0/T, \alpha, \beta, M)>0$
such that for all $a\in (0, M]$, $\lambda \in (0, T]$ and $x, y\in
\lambda^{-1}D$,
$$
p^a_{\lambda^{-1}D} (1, x, y) \,\leq  \,C_{17}\,\left(1\wedge J^a(x,
y) \right) \, \delta_{\lambda^{-1}D}(x)^{\alpha/2}.
$$
\end{thm}

\pf  Note that for every $\lambda \in (0, T]$, $\lambda^{-1}D$
satisfies the uniform exterior ball condition with radius $r_0/T$.
For $x, y\in \lambda^{-1}D$, let $z\in
\partial(\lambda^{-1} D)$
be that $|x-z|=\delta_{\lambda^{-1}D} (x)$. Let $B_z\subset
(\lambda^{-1}D)^c$ be the ball with radius $r_1:=4^{-1} \wedge
(r_0/T)$ so that $\partial B_z \cap \partial (\lambda^{-1}D)=\{z\}$.
Since, by \eqref{e:1.4}
$$
p^a_{\lambda^{-1}D} (1, x, y) \, \leq   \, p^a (1, x, y)\,\le
c\left(1 \wedge J^a(x, y) \right),
$$
it suffices to prove the theorem for $x\in \lambda^{-1}D$  with
$\delta_{\lambda^{-1}D} (x) <r_1/4$. When $\delta_{\lambda^{-1}D}
(x) <r_1/4$ and $|x-y|\geq 5r_1$, we have $\delta_{B_z^c} (y)> 2r_1$
and so by Lemma \ref{L:2.1}, there is a constant $c_1>0$ that
depends only on $(r_0/T,d, \alpha, \beta, M)$ such that
for $t\in (0, 1]$,
\begin{align}
p^a_{\lambda^{-1}D} (t, x, y)  \,\leq  \, p^a_{(\overline
B_z)^c}(t,x, y) \,\leq  \,c_1 \,  \delta_{(\overline B_z)^c}
(x)^{\alpha/2} J^a(x, y)
 = c_1 \delta_{\lambda^{-1}D}(x)^{\alpha/2} J^a(x, y)  .\label{e:2.5}
\end{align}
So it remains to show that when $\delta_{\lambda^{-1}D} (x) <r_1/4$
and $|x-y|< 5r_1$, there exists a positive constant $c_2
=c_2(r_0/T,d, \alpha, \beta, M)$ such that
 \bee\label{e:2.6}
p^a_{\lambda^{-1}D}(1, x, y)  \,\leq  \,c_2 \,
\delta_{\lambda^{-1}D}(x)^{\alpha/2}.
 \eee

Let $z_x\in \partial(\lambda^{-1} D)$  be such that
$|x-z_x|=\delta_{\lambda^{-1}D}(x)$ and $z_0\in \bR^d$ so that
$$
B(z_0,r_1)\subset (\lambda^{-1}D)^c \qquad
\hbox{and} \qquad  \partial B(z_0,r_1)
\cap \partial (\lambda^{-1}D)=\{z_x\} .
$$
Define $U:=\left\{w\in \bR^d: \, |w-z_0| \in (r_1, \,  8r_1)
\right\}$. Note that
$$
x, y \in U\cap \lambda^{-1}D \qquad \hbox{and} \qquad
\delta_U(x)=\delta_{\lambda^{-1}D}(x).
$$
By the strong Markov property and the symmetry of
$p^a_{\lambda^{-1}D}(1, x, y)$ in $x$ and $y$, we have
$$
p^a_{\lambda^{-1}D} (1, x, y) = p^a_{U\cap \lambda^{-1}D}(1, x, y) +
\E_y \left[ p^a_{\lambda^{-1}D} (1-\tau^a_{U\cap \lambda^{-1}D},
X^a_{\tau^a_{U\cap \lambda^{-1}D}}, x); \ \tau^a_{U\cap
\lambda^{-1}D}< 1  \right].
$$
By the semigroup property and \eqref{e:1.4},
\begin{eqnarray*}
p^a_{U\cap \lambda^{-1}D}(1, x, y)&=& \int_{U\cap \lambda^{-1}D}
p^a_{U\cap \lambda^{-1}D}(1/2, x, z) p^a_{U\cap \lambda^{-1}D}
(1/2, z, y) dz\\
&\leq & \| p^a(1/2, \cdot , \cdot )\|_\infty \,
\P_x \left( \tau^a_{U\cap \lambda^{-1}D} >1/2\right) \\
&\leq & c_3 \, \E_x \left[ \tau^a_{U\cap \lambda^{-1}D} \right]
\,\leq\, c_3 \, \E_x \left[ \tau^a_U  \right] \\
&\leq &  c_4\, \delta_{U}(x)^{\alpha/2}  \,= \, c_4\, \delta_{
\lambda^{-1} D}(x)^{\alpha/2}.
\end{eqnarray*}
In the last inequality, we used Lemma \ref{L:2.0}.

On the other hand, we have $X^a_{\tau^a_{U\cap \lambda^{-1}D}}\in
U^c\cap \lambda^{-1}D$ on $\{\tau^a_{U\cap \lambda^{-1}D}< 1
\}$, and so
$$
|X^a_{\tau^a_{U\cap \lambda^{-1}D}}-x|\ge 7r_1,
 \qquad \hbox{on} \quad
\left\{\tau^a_{U\cap \lambda^{-1}D}< 1 \right\}.
$$
Consequently, by \eqref{e:2.5} for $p^a_{\lambda^{-1}D}
\big(1-\tau^a_{U\cap \lambda^{-1}D},\, X^a_{\tau^a_{U\cap
\lambda^{-1}D}},\, x\big)$,
\begin{eqnarray*}
&&\E_y \left[ p^a_{\lambda^{-1}D} (1-\tau^a_{U\cap \lambda^{-1}D},
\, X^a_{\tau^a_{U\cap \lambda^{-1}D}}, \, x); \ \tau^a_{U\cap
\lambda^{-1}D}< 1 \right] \\
&&\leq   \E_y \left[ c_1\,  {\delta_{\lambda^{-1}D}(x)^{\alpha/2}}
J^a(X^a_{\tau^a_{U\cap \lambda^{-1}D}},x) ; \ \tau^a_{U\cap
\lambda^{-1}D}< 1  \right]\\
&&\leq  c_1((7r_1)^{-d-\alpha}+ M^\beta(7r_1)^{-d-\beta}) \,
\delta_{\lambda^{-1}D} (x)^{\alpha/2}\, \P_y \left( \tau^a_{U\cap
\lambda^{-1}D}< 1 \right) \\
&&\leq c_1((7r_1)^{-d-\alpha}+ M^\beta(7r_1)^{-d-\beta}) \,
\delta_{\lambda^{-1}D} (x)^{\alpha/2}.
\end{eqnarray*}
This completes the proof for \eqref{e:2.6} and hence the theorem.
\qed

\medskip

\begin{thm}\label{t:ub}
Assume that $M>0$ and that $D$ is an open set that satisfies the
uniform exterior ball condition with radius $r_0>0$. For every
$T>0$, there exists a positive constant
$C_{18}=C_{18}(T, r_0, \alpha, \beta, M)$ such that for every $a \in
[0, M]$, $t\in (0, T]$ and $x, y\in D$,
 \bee\label{e:1}
p^a_D(t, x, y) \leq C_{18} \left( 1\wedge \frac{\delta_D
(x)^{\alpha/2}}{\sqrt{t}}\right) \left( 1\wedge
\frac{\delta_D(y)^{\alpha/2}}{\sqrt{t}}\right)  \left(
t^{-d/\alpha}\wedge tJ^a(x, y)\right).
 \eee
\end{thm}

\pf Fix $T, M>0$. By Theorem \ref{ub11}, there exists a positive
constant $c_1=c_1(T, r_0, \alpha, \beta, M)$ such that for every $ t
\in (0, T]$,
 \bee\label{e:4}
p^{at^{(\alpha-\beta)/(\alpha\beta)} }_{t^{-1/\alpha}D} (1, x, y)
\leq c_1 \, \left( 1 \wedge J^{a
t^{(\alpha-\beta)/(\alpha\beta)}}(x, y) \right)
\delta_{t^{-1/\alpha}D}(x)^{\alpha/2}.
 \eee
Thus by \eqref{e:scaling}, \eqref{e:1.4} and \eqref{e:4}, for every
$t \le T$,
\begin{eqnarray*}
p^a_D(t, x, y) &=& t^{-d/\alpha} p^{at^{(\alpha - \beta)/
(\alpha\beta)}}_{t^{-1/\alpha} D} (1, t^{-1/\alpha}x,
t^{-1/\alpha} y)\\
&\leq & c_1 \, t^{-d/\alpha}  \left( 1 \wedge J^{a
t^{(\alpha-\beta)/(\alpha\beta)}}(t^{-1/\alpha}x, t^{-1/\alpha}
y)\right) \delta_{t^{-1/\alpha} D}(t^{-1/\alpha}
x)^{\alpha/2} \\
&=& c_1  \left( t^{-d/\alpha} \wedge tJ^a(x, y) \right)
\frac{\delta_{D}( x)^{\alpha/2}}{\sqrt{t}} \\
&\leq & c_2 \, p^a(t, x, y)\, \frac{\delta_{D}(
x)^{\alpha/2}}{\sqrt{t}}.
\end{eqnarray*}
By  symmetry, the  above  inequality  holds with the roles of $x$
and $y$ interchanged. Using the semigroup property for $t\leq T$,
\begin{eqnarray*}
p^a_D(t, x, y) &=& \int_D p^a_D(t/2, x, z) p^a_D (t/2, z, y) dz\\
&\leq & c_3\,  \frac{\delta_{D}( x)^{\alpha/2} \delta_{D}(
y)^{\alpha/2}}{t } \int_D p^a(t/2, x, z) p^a (t/2, z, y) dz \\
&\leq & c _3\, \frac{\delta_{D}( x)^{\alpha/2} \delta_{D}(
y)^{\alpha/2}}{t } p^a(t, x, y) .
\end{eqnarray*}
This proves the upper bound \eqref{e:1} by noting that
$$
(1\wedge a)(1\wedge b) = \min\{1, a, b, ab\} \qquad \hbox{for } a,
b>0.
$$
\qed

\section{Lower bound estimate}

\begin{lemma}\label{L:4.2}
For any positive constants  $\Lambda$, $\kappa$ and $b$, there exists
$C_{19}=C_{19}(\Lambda, \kappa, b, \alpha, \beta, M)>0$ such that
for every $z \in \bR^d$, $\lambda \in (0, \Lambda]$ and $a\in (0,
M]$,
$$
\inf_{y\in\bR^d \atop |y -z| \le \kappa \lambda^{1/\alpha}} \P_y \left(
\tau^a_{B(z, 2\kappa \lambda^{1/\alpha} )} > b\lambda \right)\, \ge\,
C_{19}.
$$
\end{lemma}

\pf By \cite[Proposition 4.9]{CK2}, there exists $\eps=\eps(\Lambda,
\kappa, \alpha, \beta)>0$ such that for every
$\lambda \in (0, \Lambda]$,
$$
\inf_{y\in \bR^d } \P_y \left( \tau^1_{B(y, \kappa \lambda^{1/\alpha}/2
)}
> \eps \lambda \right) \ge \frac12.
$$
Suppose $b>\eps$ then by the parabolic Harnack principle in
\cite[Proposition 4.12]{CK2}
$$
c_1\,p ^1_{B(y, \kappa \lambda^{1/\alpha} )}(\eps \lambda ,y,w) \, \le \,
p ^1_{B(y, \kappa \lambda^{1/\alpha} )}(b \lambda ,y,w)\qquad \hbox{for }
w \in B(y, \kappa \lambda^{1/\alpha}/2 ),
$$
where the constant $c_1=c_1(\kappa, b, \alpha, \beta, \Lambda)>0$ is
independent of $y\in \bR^d$, $\lambda\in (0, \Lambda]$. Thus
\begin{align}
\P_y \left( \tau^1_{B(y, \kappa \lambda^{1/\alpha} )} > b\lambda \right)
&=\int_{B(y, \kappa \lambda^{1/\alpha} )}
p ^1_{B(y, \kappa \lambda^{1/\alpha} )}(b\lambda ,y,w) dw\nn\\
&\ge \int_{B(y, \kappa \lambda^{1/\alpha}/2 )} p ^1_{B(y, \kappa
\lambda^{1/\alpha} )}(b \lambda ,y,w) dw\nn\\
&\ge c_1 \int_{B(y, \kappa  \lambda^{1/\alpha}/2 )} p^1_{B(y, \kappa
\lambda^{1/\alpha}/2 )}(\eps \lambda ,y,w) dw \, \ge\, c_1/2. \label{e:dds}
\end{align}
For the general case, by \eqref{e:scaling} and \eqref{e:dds},
\begin{align*}
&\inf_{y\in\bR^d \atop |y -z| \le \kappa \lambda^{1/\alpha}} \P_y \left(
\tau^a_{B(z, 2\kappa \lambda^{1/\alpha} )} > b\lambda \right)  \\
&\ge
\P_0 \left( \tau^a_{B(0, \kappa \lambda^{1/\alpha} )} > b\lambda \right) \\
&= \int_{B(0, \kappa \lambda^{1/\alpha} )}
p ^a_{B(0, \kappa \lambda^{1/\alpha} )}(b\lambda ,0,w) dw\\
&= \int_{B(0, \kappa \lambda^{1/\alpha} a^{-\frac{\beta}{\alpha-\beta}}
)} p ^1_{B(0, \kappa \lambda^{1/\alpha}a^{-\frac{\beta}{\alpha-\beta}})}
(a^{\frac{\alpha\beta}{\alpha-\beta}} b\lambda ,0,z) dz\\
&=\P_0 \left( \tau^1({B(0, \kappa \lambda^{1/\alpha} a^{-\frac{\beta}
{\alpha-\beta}})}) >a^{\frac{\alpha\beta}{\alpha-\beta}} b\lambda
\right)\\
&\ge \P_0 \left( \tau^1(B(0, \kappa \lambda^{1/\alpha} M^{-\frac{\beta}
{\alpha-\beta}})) >M^{\frac{\alpha\beta}{\alpha-\beta}} b\lambda
\right)
\ge  c_2(\Lambda, \kappa, b, \alpha, \beta, M)
>0.
\end{align*}
This proves the lemma.   \qed

Recall that $\psi^a$ is defined in \eqref{e:psi}.

\begin{prop}\label{step1}
Suppose that $M, T>0$ and $(t, x, y)\in (0, T]\times D\times D$ with
$\delta_D(x) \ge t^{1/\alpha} \geq 2|x-y| \psi^a
(|x-y|)^{-1/(d+\alpha)} $. Then there exists a positive constant
$C_{20}=C_{20}(M, \alpha, \beta,  T)$ such that for all $a\in (0,
M]$
 \bee\label{e:lb1}
p^a_D(t,x,y) \,\ge\,C_{20}\, t^{-d/\alpha}.
 \eee
\end{prop}

\pf Let $t\in (0, T]$ and $x, y \in D$ with $\delta_D(x) \ge
t^{1/\alpha} \geq 2|x-y| \psi^a(|x-y|)^{-1/(d+\alpha)} $. By the
parabolic Harnack principle in \cite[Proposition 4.12]{CK2} and the
scaling property, there exists $c_1=c_1(M, \alpha, \beta , T)>0$
such that for all $a\in (0, M]$,
$$
p^a_D(t/2,x,w) \, \le  \, c_1\, p^a_D(t,x,y)\qquad \hbox{for } w \in
B(x, 2t^{1/\alpha}/3) .
$$
This together with Lemma \ref{L:4.2} yields that
\begin{eqnarray*}
p^a_D(t, x, y) &\geq & \frac{1}{c_1 \, | B(x,t^{1/\alpha}/2)|}
\int_{B(x,t^{1/\alpha}/2)} p^a_D(t/2, x, w)dw\\
&\geq & c_2 t^{-d/\alpha} \, \int_{B(x,t^{1/\alpha}/2)}
p^a_{B(x,t^{1/\alpha}/2)} (t/2, x, w)dw \\
&=& c_2 t^{-d/\alpha} \, \P_x \left( \tau^a_{B(x,t^{1/\alpha}/2)} >
t/2\right) \ \geq \ c_3 \, t^{-d/\alpha},
\end{eqnarray*}
where $c_i=c_i(T, \alpha, \beta,  M)>0$ for $i=2, 3$. \qed

\begin{lemma}\label{l:st3_3}
Suppose that
$M, T>0$,  $D$ is an open subset of $\bR^d$ and $(t, x, y)\in
(0,T]\times D\times D$ with $\min\left\{ \delta_D(x), \, \delta_D
(y)\right\} \ge t^{1/\alpha}$ and $|x-y|^\alpha \ge 2^{-\alpha} t
\psi^a (|x-y|)^{\alpha/(d+\alpha)}$. Then  there exists a constant
$C_{21}=C_{21}(\alpha, \beta, T, M)>0$ such that for $a\in (0, M]$
$$
\P_x \left( X^{a,D}_t \in B \big( y, \,  2^{-1} t^{1/\alpha} \big)
\right) \,\ge \, C_{21}\, t^{d/\alpha+1} J^a(x,y).
$$
\end{lemma}

\pf
For $t\in (0, T]$, it follows from Lemma \ref{L:4.2} that, starting
at $z\in B(y, \, 4^{-1} t^{1/\alpha})$, with probability at least
$c_1=c_1(\alpha,
\beta, T, M)>0$, for any $a\in (0, M]$,  the process $X^a$ does not
move more than $6^{-1}t^{1/\alpha}$ by time $t$. Thus, it suffices
to show that there exists a constant $c_2=c_2(\alpha, \beta, T,
M)>0$ such that
\begin{equation}\label{eq:molow}
\P_x \left(X^{a,D} \hbox{ hits the ball } B(y, \,
4^{-1}t^{1/\alpha})\mbox{ by time } t \right)\, \ge\, c_2\,
{t^{d/\alpha+1}}{ J^a(x,y)}
\end{equation}
for all $a\in (0, M]$, $t \in (0, T] $ and  $|x-y|^\alpha \ge
2^{-\alpha} t \psi^a (|x-y|)^{\alpha/(d+\alpha)}$.

Let $B_x:=B(x, \, 6^{-1}t^{1/\alpha})$, $B_y:=B(y, \, 6^{-1}
t^{1/\alpha})$ and $\tau^a_x:=\tau^a_{B_x}$. It follows from Lemma
\ref{L:4.2} that there exists $c_3=c_3(\alpha, \beta, T, M)>0$ such
that for $a\in (0, M]$ and $t\in (0, T]$,
\begin{equation}\label{eq:lowtau}
\E_x \left[ t \wedge\tau^a_{x} \right] \,\ge\,
{t} \P_x\left(\tau^a_{x} \ge t\right)
\,\ge\,  c_3\,t .
\end{equation}
By the L\'evy system in \eqref{e:levy},
\begin{eqnarray}
&&\P_x \left(X^{a,D} \mbox{ hits the ball } B(y, \,
4^{-1} t^{1/\alpha} )\mbox{ by time }   t \right) \nonumber\\
&\ge & \P_x(X^a_{t\wedge \tau^a_x}\in B(y,\, 4^{-1} t^{1/\alpha})
\hbox{ and } t \wedge \tau^a_x \hbox{ is a jumping time })\nonumber \\
&\geq & \E_x \left[\int_0^{t\wedge \tau^a_x} \int_{B_y} J^a(X^a_s,
u) duds \right] \label{e:dew}.
\end{eqnarray}
Note that
$$
|x-y| \ge 2^{-1} t^{1/\alpha} \psi^a (|x-y|)^{1/(d+\alpha)}\ge
2^{-1} t^{1/\alpha}.
$$
Moreover, if $s < \tau^a_{x}$ and $u \in B_y$,
$$
|X^a_s-u|\le|x-y|+|x-X^a_s|+|y-u| \le 2|x-y|.
$$
Thus from \eqref{e:dew} we get that for any $a\in (0, M]$ and $t\in (0, T]$,
\begin{eqnarray*}
&&\P_x \left(X^{a,D} \mbox{ hits the ball } B(y, \,
4^{-1} t^{1/\alpha} )\mbox{ by time }   t \right)\\
&\ge & \E_x \left[ t\wedge \tau^a_x \right] \int_{B_y}
j^a(2|x-y|)\,  du\\
&\ge&   c_4\,  t\, |B_y| \,  j^a(2|x-y|)  \,\ge\, c_5\,
{t^{d/\alpha+1}}{ j^a(2|x-y|)} \,\ge\, c_5\, 2^{-d-\alpha}
{t^{d/\alpha+1}}{ j^a(|x-y|)}
\end{eqnarray*}
for some positive constants $c_i=c_i(\alpha, \beta, T, M)$, $i=4,
5$. Here in the second inequality, (\ref{eq:lowtau}) is used. \qed

\begin{prop}\label{step3}
Suppose that $T>0$, $M>0$, $D$ is an open subset of $\bR^d$
and $(t, x, y)\in (0, T]\times D\times D$ with
$\min\left\{ \delta_D(x), \, \delta_D (y)\right\} \ge
(t/2)^{1/\alpha}$ and $|x-y|^\alpha \ge 2^{-\alpha-1} t \psi^a
(|x-y|)^{\alpha/(d+\alpha)}$. Then there exists a constant
$C_{22}=C_{22}(\alpha, \beta, T, M)>0$ such that for any $a\in (0,
M]$,
\begin{equation}\label{e:lu3}
p^a_D(t, x, y)\,\ge  \,C_{22} \,{ t}{J^a(x,y)}.
\end{equation}
\end{prop}

\pf By the semigroup property,  Proposition \ref{step1} and Lemma
\ref{l:st3_3}, there exist positive constants $c_1=c_1(\alpha,
\beta, T, M)$ and $c_2=c_2(\alpha, \beta, T, M)$ such that for any
$t\in (0, T]$ and $a\in (0, M]$
\begin{eqnarray*}
p^a_D(t, x, y) &=& \int_{D} p^a_{D}(t/2, x, z) p^a_{D}(t/2, z,
y)dz\\
&\ge& \int_{B(y, \, 2^{-1} (t/2)^{1/\alpha})} p^a_{D}(t/2, x, z)
p^a_{D}(t/2, z, y) dz\\
&\ge& c_1 t^{-d/\alpha} \bP_x \left( X^{a,D}_{t/2} \in B(y, 2^{-1}
(t/2)^{1/\alpha}) \right) \\
& \ge & c_2\,{ t}{J^a(x,y)}.
\end{eqnarray*}
 \qed
\medskip

In the rest of this section, we assume that $D$ is an open set in
$\bR^d$ satisfying the uniform interior ball condition with radius
$r_0>0$ in the following sense: For every $x\in D$ with $\delta_D
(x)<r_0$, there is $z_x\in \partial D$ so that $|x-z_x|=\delta_D(x)$
and $B(x_0, r_0)\subset D$ for $x_0:= z_x+r_0 (x-z_x)/|x-z_x|$.
Clearly, a (uniform) $C^{1,1}$ open set satisfies the uniform
interior ball condition.

The goal of this section is to prove the following lower bound for
the heat kernel $p^a_D(t, x, y)$.

\begin{thm}\label{t:lb}
For any $M>0$ and  $T>0$, there exists  positive constant
$C_{23}=C_{23}(\alpha, \beta, T, M, r_0)$ such that  for all $a\in
(0, M]$ and $(t, x, y)\in (0, T]\times D\times D$,
$$
p^a_D(t, x, y) \ge C_{23} \left( 1\wedge
\frac{\delta_D(x)^{\alpha/2}}{\sqrt{t}}\right) \left( 1\wedge
\frac{\delta_D(y)^{\alpha/2}}{\sqrt{t}}\right) \left( t^{-d/\alpha}
\wedge tJ^a(x, y) \right).
$$
\end{thm}

To prove this result, we will first prove a lower bound estimates on
the Green function of $X^{a, U}$
$$
G^a_U(x, y):=\int_0^\infty p^a_U(t, x, y)dt
$$
when $U$ is a bounded $C^{1, 1}$ open set. The tool we use to
establish the Green function lower bound is a subordinate killed
$\alpha$-stable process in $U$. We first introduce this subordinate
killed process first.

Assume that $U$ is a bounded $C^{1, 1}$ open set in $\R^d$ and $R_1$
the radius in the uniform interior and exterior ball conditions.
Then it follows from \cite[Theorem 1.1]{CKS} that the killed
$\alpha$-stable process $X^U$ on $U$ has a density $p_U(t, x, y)$
satisfying the following condition: for any $T>0$ there exist
positive constants $c_2>c_1$ depending only on $\alpha, T, R_1$ and
$d$ such that for any $(t, x, y)\in (0, T]\times U\times U$,
\begin{eqnarray}
p_U(t, x, y)&\ge & c_1\left( 1\wedge
\frac{\delta_U(x)^{\alpha/2}}{\sqrt{t}}\right) \left( 1\wedge
\frac{\delta_U(y)^{\alpha/2}}{\sqrt{t}}\right) \left( t^{-d/\alpha}
\wedge \frac{t}{|x-y|^{d+\alpha}}\right),\label{e:kstlb}\\
p_U(t, x, y)&\le & c_2\left( 1\wedge
\frac{\delta_U(x)^{\alpha/2}}{\sqrt{t}}\right) \left( 1\wedge
\frac{\delta_U(y)^{\alpha/2}}{\sqrt{t}}\right) \left( t^{-d/\alpha}
\wedge \frac{t}{|x-y|^{d+\alpha}}\right).\label{e:kstub}
\end{eqnarray}
Let $\{T^a_t: t\ge 0\}$ be a subordinator, independent of $X^a$,
with Laplace exponent
$$
\phi^a(\lambda)=\lambda + a^\beta \lambda^{\beta/\alpha}.
$$
Then the process $\{Z^{a, U}_t: t\ge 0\}$ defined by $Z^{a,
U}_t=X^U_{T^a_t}$ is called a subordinate killed stable process in
$U$. Since $\phi^a$ is a complete Bernstein function, the
subordinate $T^a$ has a decreasing potential density $u^a(x)$. In
fact $u^a(x)$ is completely monotone. (See \cite{RSV, SV06} for the
details.) Then it follows from \cite{SV06} that the Green function
$R^a_U(x, y)$ of $Z^{a, U}$ is given by
\begin{equation}\label{gfn4sks}
R^a_U(x, y)=\int^{\infty}_0p_U(t, x, y)u^a(t)dt.
\end{equation}
It follows from \cite{SV08} that the Green function $G^a_U$ of
$X^{a, U}$ and the Green function $R^a_U$ of $Z^{a, U}$ satisfy the
following relation:
\begin{equation}\label{gfcrln}
R^a_U(x, y)\le G^a_U(x, y) \qquad (x, y)\in U\times U.
\end{equation}
So we can get a lower bound on $G^a_U(x, y)$ be establishing a lower
bound on $R^a_U(x, y)$. The following result gives sharp two-sided
estimates on $R^a_U(x, y)$ and the idea of the proof is similar to
that of \cite{SV}.

\begin{thm}\label{T:gfcesks}
Suppose that $M>0$ and $U$ is a bounded $C^{1, 1}$ open set in $U$.
There exist positive constant $C_{25}>C_{24}$ depending only on
$(\alpha, \beta, d, R_1, M, \text{diam}(U))$ such that for all $a\in
(0, M]$,
$$
R^a_U(x, y)\ge C_{24} \begin{cases}
 \left(1\wedge \frac{ \delta_U(x)^{\alpha/2} \delta_U(y)^{\alpha/2}}{
|x-y|^{\alpha}}
\right) |x-y|^{\alpha -d}   \qquad &\hbox{when } d>\alpha ,  \\
\log \left( 1+ \frac{  \delta_U(x)^{\alpha/2} \delta_U
(y)^{\alpha/2}}{ |x-y|^{\alpha}
}\right)   &\hbox{when } d=1=\alpha , \\
\big( \delta_U(x)  \delta_U (y)\big)^{(\alpha-1)/2} \wedge \frac{
\delta_U(x)^{\alpha/2} \delta_U (y)^{\alpha/2}}{ |x-y|}  &\hbox{when
} d=1<\alpha ,
\end{cases}
$$
and
$$
R^a_U(x, y)\le C_{25} \begin{cases}
 \left(1\wedge \frac{  \delta_U(x)^{\alpha/2}
\delta_U(y)^{\alpha/2}}{ |x-y|^{\alpha}}
\right) |x-y|^{\alpha -d}   \qquad &\hbox{when } d>\alpha ,  \\
\log \left( 1+ \frac{  \delta_U(x)^{\alpha/2} \delta_U
(y)^{\alpha/2}}{ |x-y|^{\alpha}
}\right)   &\hbox{when } d=1=\alpha , \\
\big( \delta_U(x)  \delta_U (y)\big)^{(\alpha-1)/2} \wedge \frac{
\delta_U(x)^{\alpha/2} \delta_U (y)^{\alpha/2}}{ |x-y|}  &\hbox{when
} d=1<\alpha .
\end{cases}
$$
\end{thm}

\pf Since the drift coefficient of $T^a$ is 1, we know that
$u^a(t)\le 1$ for all $t>0$. Now the upper bound on $R^a_U$ follows
immediately from \eqref{gfn4sks} and \cite[Corollary 1.2]{CKS}.
 Thus we only need to prove the lower bound.

By using a scaling argument, one can easily check that
\begin{equation}\label{e:relpd}
u^a(t)=u^1(a^{\frac{\alpha}{\alpha-\beta}}t) \qquad t>0.
\end{equation}
Let $T={\rm diam}(U)$. Since $u^1(t)$ is a completely monotone
function with $u^1(0+)=1$, by \eqref{e:relpd},
\begin{align}
u^a(t)\ge  u^1(M^{\frac{\alpha}{\alpha-\beta}}T)  \qquad \hbox{for every }
t\in (0, T] \hbox{ and } a\in (0, M].
\label{e:stlb4pd}
\end{align}
Using \eqref{e:stlb4pd}, \eqref{gfn4sks} and \cite[(4.2)]{CKS} we
get that
\begin{eqnarray}
&& \int_0^T\left( 1\wedge
\frac{\delta_U(x)^{\alpha/2}}{\sqrt{t}}\right) \left( 1\wedge
\frac{\delta_U(y)^{\alpha/2}}{\sqrt{t}}\right) \left( t^{-d/\alpha}
\wedge \frac{t}{|x-y|^{d+\alpha}}\right) u^a(t)dt \nonumber\\
&\ge&\frac{u^1(M^{\frac{\alpha}{\alpha-\beta}}T)}{|x-y|^{d-\alpha}}
\int_{\frac{|x-y|^\alpha}{T}}^\infty \left(u^{\frac{d}{\alpha}-2}
\wedge u^{-3} \right)  \left(1\wedge \frac{ {\sqrt u}
\delta_U(x)^{\alpha/2}}{ |x-y|^{\alpha/2} }\right) \left(1\wedge
\frac{ {\sqrt u} \delta_U(y)^{\alpha/2}}{ |x-y|^{\alpha/2} }\right)
du.          \label{ID:5.1}
\end{eqnarray}
Now we can follow the proof of \cite[Corollary 1.2]{CKS} to get the
desired lower bound. In fact, when $d>\alpha$, the desired lower
bound follows from \eqref{ID:5.1} and \cite[(4.3) and (4.7)]{CKS}.
Let
$$
u_0:=  \frac{ \delta_U(x)^{\alpha/2} \delta_U(y)^{\alpha/2}}{
|x-y|^{\alpha} }.
$$
When $d=\alpha=1$, by \eqref{ID:5.1} and \cite[(4.3) and (4.9)]{CKS},
\begin{eqnarray*}
R^a_U(x, y)&\ge &u^1(M^{\frac{\alpha}{\alpha-\beta}}T)\int_0^Tp_U(t, x, y) dt \\
&\ge& c_1 \left(1\wedge \frac{  \delta_U(x)^{\alpha/2}}{
|x-y|^{\alpha/2}   }\right) \left( 1\wedge \frac{
\delta_U(y)^{\alpha/2}}{ |x-y|^{\alpha/2}} \right) +c_1 \log (u_0 \vee
1) + c_1u_0 \left( (1/u_0)\wedge 1 -
\frac{|x-y|^\alpha}{T} \right)  \\
&\ge& c_2(1\wedge u_0) + c_2\log (u_0 \vee 1) +c_2 u_0 \left( (1/u_0)\wedge
1 - \frac{|x-y|^\alpha}{T} \right)  \\
&\ge& c_3(1\wedge u_0) + c_3\log (u_0 \vee 1)\,\ge \, c_4 \log
\left( 1+ \frac{  \delta_U(x)^{\alpha/2} \delta_U
(y)^{\alpha/2}}{ |x-y|^{\alpha} }\right) .
\end{eqnarray*}
Lastly, in the case $d=1<\alpha<2$. By \eqref{ID:5.1}, \cite[(4.3)
and (4.7)]{CKS} and  the first display in
part (iii) of the proof of \cite[Corollary1.2]{CKS}, we have
\begin{eqnarray*}
R^a_U(x, y)&\ge&u^1(M^{\frac{\alpha}{\alpha-\beta}}T)
\int_T^\infty p_U(t, x, y)dt \\
&\ge &
c_5\frac1{|x-y|^{1-\alpha}} \left(1\wedge u_0 \right)     \\
&&+ c_5\frac{1}{|x-y|^{1-\alpha}}\left( \left( (u_0\vee
1)^{1-(1/\alpha)}-1\right) + c_5  u_0 \left( (u_0\vee 1)^{-1/\alpha}
-\left(\frac{|x-y|^\alpha}T\right)^{1/\alpha}\right)
\right) \\
&\ge&
c_6\frac{1}{|x-y|^{1-\alpha}}\left(  u_0\wedge u_0^{1-(1/\alpha)}\right)\\
&=&
c_6\left(\left( \delta_U(x) \delta_U (y)\right)^{(\alpha-1)/2} \,
\wedge \, \frac{  \delta_U(x)^{\alpha/2} \delta_U (y)^{\alpha/2}}{
|x-y| } \right).
\end{eqnarray*}
 \qed

By integrating the lower bound in Theorem \ref{T:gfcesks} with
respect to $y$ and applying \eqref{gfcrln}, we obtain the following
lower bound on $\E_x[\tau^a_U]$
\begin{cor}\label{C:taul}
Suppose that $M>0$ and $U$ is a bounded $C^{1, 1}$ open set in $U$.
Then there exists a constant
$C_{26}=C_{26}(\alpha, \beta, d, M, R_1, \text{diam}(U))>0$ such
that for every $a \in (0, M]$ and $x\in U$,
$$
\E_x[\tau^a_U] \ge C_{26}  \delta_U(x)^{\alpha/2} .\
$$
\end{cor}

We will first establish Theorem \ref{t:lb} for small $T$, that is,
we will first assume that
 \bee\label{tsmall}
 t\le T_0:= \left(\frac{r_0}{16}\right)^{\alpha}.
 \eee
By integrating \eqref{e:scaling} with respect to $t$ and $y$, we
have that for every open set $U$, $\lambda>0$ and $x\in U$,
\begin{equation}\label{e:scaling2}
\E_x[\tau^a_U]= \int_{U} G^{a}_U (  x, z)  dz=
\lambda^{\alpha}\int_{\lambda^{-1}U}
G^{a\lambda^{(\alpha-\beta)/\beta}}_{\lambda^{-1}U} (\lambda^{-1} x,
y)dy = \lambda^{\alpha}
\E_{\lambda^{-1}x}
\left[\tau^{a\lambda^{(\alpha-\beta)/\beta}}({\lambda^{-1}U})\right].
\end{equation}

\medskip

\begin{lemma}\label{l:st2_3}
Suppose that  $M>0$, $\kappa \in (0, 1)$ and that $(t, x)\in (0,
T_0]\times D$ with $\delta_D(x) \leq 3 t^{1/\alpha} < r_0/4$. Let
$z_x\in \partial D$ be such that $|z_x-x|=\delta_D (x)$ and define
${\bf n}(z_x):=(x-z_x)/|x-z_x|$. Put $x_1=z_x+ 3t^{1/\alpha}{\bf
n}(z_x)$ and $B=B(x_1, 3t^{1/\alpha})$. Suppose that $x_0$ is a
point on the line segment connecting $z_x$ and $z_x+
6t^{1/\alpha}{\bf n}(z_x)$ such that $B(x_0, 2 \kappa
t^{1/\alpha})\subset B \setminus \{x\}$. Then for any $b>0$, there
exists a constant
$C_{27}=C_{27}(\kappa, \alpha, \beta, r_0, b, M)>0$ such that for
all $a\in (0, M]$
\begin{equation}\label{e:4.5}
\P_x \left( X^{a,D}_{bt} \in B(x_0, \kappa t^{1/\alpha}) \right) \,\ge
\, C_{27}\, t^{-1/2} \delta_D(x)^{\alpha/2}.
\end{equation}
\end{lemma}

\pf Let $0< \kappa_1 \le \kappa$ and assume first that
$2^{-4}\kappa_1 t^{1/\alpha} < \delta_D(x)\leq 3t^{1/\alpha}$.
Repeating the proof of Lemma \ref{l:st3_3}, we get that, in this
case, there exists a constant $c_1=c_1(\alpha, \beta, \kappa_1, M,
r_0, b)>0$ such that for all $a\in (0, M]$
$$
\P_x \left( X^{a,D}_{bt} \in B(x_0, \kappa_1 t^{1/\alpha}) \right) \,\ge
\,c_1t^{d/\alpha+1}J^a(x, x_0)
\,\ge
\,c_1 \sA(d, -\alpha) t^{d/\alpha+1}|x- x_0|^{-d-\alpha}
$$
for all $t\le T_0$. Using the fact that $|x-x_0|\in [2\kappa
t^{1/\alpha}, 6t^{1/\alpha}]$ we get that for all $a\in (0, M]$,
\begin{equation}\label{e:case1}
\P_x \left( X^{a,D}_{bt} \in B(x_0, \kappa_1 t^{1/\alpha}) \right) \,\ge
\,c_2>0
\end{equation}
for some constant $c_2=c_2(\alpha, \beta, \kappa_1, M, r_0, b)$. By
taking $\kappa_1=\kappa$, this shows that \eqref{e:4.5} holds for
all $b>0$ in the case when $2^{-4}\kappa_1 t^{1/\alpha} <
\delta_D(x)\leq 3t^{1/\alpha}$.

So it suffices to consider the case that $\delta_D(x) \leq
2^{-4}\kappa t^{1/\alpha}$. We now show that there is some $b_0>1$
so that \eqref{e:4.5} holds for every $b\geq b_0$ and $\delta_D(x)
\leq 2^{-4}\kappa t^{1/\alpha}$. For simplicity, we assume without
loss of generality that $x_0=0$ and let $\wh B:=B(0, \kappa
t^{1/\alpha})$.  Let $x_2=z_x+4^{-1}\kappa {\bf n}(z_x)t^{1/\alpha}$
and  $B_2:= B(x_2, 4^{-1}\kappa t^{1/\alpha})$. Observe that since
$B(0, 2\kappa t^{1/\alpha}) \subset B \setminus \{x\}$,
\begin{equation}\label{e:4.7}
\kappa/2 t^{1/\alpha}\leq |y-z| \leq 6 t^{1/\alpha}  \qquad
\hbox{for } y \in B_2 \hbox{ and } z\in B(0, \kappa t^{1/\alpha}).
\end{equation}
By the strong Markov property of $X^a$ at the first exit time
$\tau^a_{B_2}$ from $B_2$ and Lemma \ref{L:4.2},
\begin{eqnarray}
&& \P_{x} \left(X^a_{bt} \in B(0, \kappa t^{1/\alpha}) \right) \nonumber \\
&\geq & \P_{x}\left( \tau^a_{B_2}<bt, X^a_{\tau^a_{B_2}} \in B(0,
2^{-1}\kappa t^{1/\alpha}) \hbox{ and }
|X^a_s-X^a_{\tau_{B_2}}|<\kappa/2 \hbox{ for } s\in
[\tau^a_{B_2}, \tau^a_{B_2}+bt^{1/\alpha}]\right) \nonumber \\
&\geq& c_3 \, \P_{x}\left( \tau^a_{B_2}<bt \hbox{ and }
X^a_{\tau^a_{B_2}}\in B(0, 2^{-1}\kappa t^{1/\alpha})
\right)\label{e:4.8}.
\end{eqnarray}

It follows from  \eqref{e:scaling2} and Corollary
\ref{C:taul} that
\begin{eqnarray}
\P_{x} \left(  X^a_{\tau^a_{B_2}}\in B(0, 2^{-1}\kappa
t^{1/\alpha})\right)
&=&\int_{B(0, 2^{-1}\kappa t^{1/\alpha})}\int_{B_2}G^a_{B_2}(x,
y)J^a(y, z)dydz\nonumber\\
&\ge & \sA (d, -\alpha)\int_{B(0, 2^{-1}\kappa t^{1/\alpha})}\int_{B_2}G^a_{B_2}(x,
y)\frac{dydz}{|y-z|^{d+\alpha}}\nonumber\\
&\ge&\frac{c_4}{t}\E_x\left[\tau^a_{B_2}\right] \nn\\
&=&c_4
\E_{x/t^{1/\alpha}}
\left[\tau^{at^{(\alpha-\beta)/{\alpha\beta}}}({B(x_2/t^{1/\alpha},
4^{-1}\kappa)}) \right]\nonumber\\
&\ge&c_5\left(\frac{\delta_{D}(x)}{t^{1/\alpha}}\right)^{\alpha/2}\label{e:4.9}
\end{eqnarray}
for some positive constants $c_4, c_5$ depending only on $\alpha,
\beta, r_0, \kappa$ and $M$. Note that, by \eqref{e:scaling}
$$
\int_{B(x_2, 4^{-1}\kappa t^{1/\alpha})} p^a_{ B(x_2, 4^{-1}\kappa
t^{1/\alpha})} (bt, x, z)dz = \int_{B(t^{-1/\alpha} x_2,
4^{-1}\kappa )} p^{at^{(\alpha-\beta)/\alpha\beta}}_{
B(t^{-1/\alpha} x_2, 4^{-1}\kappa)} (b, t^{-1/\alpha}x, w)dw.
$$
Since $at^{(\alpha-\beta)/\alpha\beta} \le
MT_0^{(\alpha-\beta)/\alpha\beta}$,  by applying Theorem \ref{t:ub}
to the right hand side of the above display,  we get
\begin{eqnarray}
\P_{x}(  \tau^a_{B_2}\geq b t ) &\leq&
b^{-d/\alpha} \int_{B(t^{-1/\alpha} x_2, 4^{-1}\kappa )}
\frac{\delta_{ B(t^{-1/\alpha} x_2, 4^{-1}\kappa)}
(t^{-1/\alpha}x)^{\alpha/2}}{ {\sqrt b}} dw \nn\\
&\leq&c_6  \ b^{-d/\alpha-1/2} \  \delta_{ t^{-1/\alpha} D}
(t^{-1/\alpha}x)^{\alpha/2} = c_6 \ b^{-d/\alpha-1/2} \,
\left(\frac{\delta_{D}(x)}{t^{1/\alpha}}
\right)^{\alpha/2},\label{e:caseaa}
\end{eqnarray}
for some positive constant $c_6$ depending only on
$\alpha, \beta, r_0, \kappa$  and  $M$.
Define
$$
b_0:= \left(\frac{2c_6}{c_5}\right)^{\frac{2\alpha}{2d+\alpha}}.
$$
We have by \eqref{e:4.8}--\eqref{e:caseaa} that for $b\geq b_0$,
\begin{eqnarray}
\P_x ( X^a_{bt}\in \wh B ) &\geq& c_3 \, \left( \P_{x}  \big(
X^a_{\tau^a_{B_2}}\in B(0, 2^{-1}\kappa t^{1/\alpha} ) \big)
-  \P_{x} \left(  \tau^a_{B_2}\geq bt  \right) \right) \nonumber \\
&\geq&  c_3 \, (c_5/2) \, \left(\frac{\delta_{D}(x)}{t^{1/\alpha}}
\right)^{\alpha/2}.\label{e:case2}
\end{eqnarray}
\eqref{e:case1} and \eqref{e:case2} show that \eqref{e:4.5} holds
for every $b\geq b_0$ and for every $x\in D$ with $\delta_D(x) \leq
3 t^{1/\alpha}$.

Now we deal with the case $0<b<b_0$ and $\delta_D(x) \leq
2^{-4}\kappa t^{1/\alpha}$. If $\delta_D(x) \leq 3
(bt/b_0)^{1/\alpha}$, we have from \eqref{e:4.5} for the case of
$b=b_0$ that
\begin{eqnarray*}
\P_x \left( X^a_{bt}\in B(x_0,  \,   \kappa t^{1/\alpha}) \right)
&\ge & \P_x \left( X^a_{b_0 (bt/b_0)}\in B(x_0,   \,
\kappa  (bt/b_0)^{1/\alpha}) \right) \\
&\geq&  c_7 \, \left(\frac{\delta_D(x)}{(bt/b_0)^{1/\alpha}}
\right)^{\alpha/2}= c_8 \,
\left(\frac{\delta_{D}(x)}{t^{1/\alpha}}\right)^{\alpha/2}.
\end{eqnarray*}
If $3 (bt/b_0)^{1/\alpha} < \delta_D(x) \leq 2^{-4}\kappa
t^{1/\alpha}$ (in this case $\kappa > 3 \cdot 2^4
(b/b_0)^{1/\alpha}$), we get \eqref{e:4.5} from \eqref{e:case1} by
taking $\kappa_1=(b/b_0)^{1/\alpha}$. The proof of the lemma is now
complete. \qed

\medskip

\begin{prop}\label{step2}
Suppose that $M>0$ and $(t, x, y)\in (0, T_0]\times D\times D$ with
$|x-y| \leq t^{1/\alpha}\psi^a(|x-y|)^{1/(d+\alpha)}$, $\delta_D
(x)\leq 2 t^{1/\alpha}$ and $\delta_D(y)\le r_0/5$. Then there
exists a constant $C_{28}=C_{28}(\alpha, \beta, M, r_0)>0$  such
that for all $a\in (0, M]$,
\begin{equation}\label{e:lu2}
p^a_D(t, x, y)\,\ge  \,C_{28} \,t^{-d/\alpha-1}\,
\delta_D(x)^{\alpha/2} \delta_D(y)^{\alpha/2} .
\end{equation}
\end{prop}

\pf Under the assumptions of the proposition, there are points $z_x,
\, z_y\in \partial D$ and $x_0, \, y_0 \in D$ such that
$\delta_D(x)=|x-z_x|$,  $\delta_D (y)=|y-z_y|$, $\partial B(x_0,
4t^{1/\alpha}) \cap \partial D = \{z_x\}$ and $\partial B(y_0,
4t^{1/\alpha})\cap \partial D= \{z_y\}$. Observe that
$$
\delta_D(x_0)=\delta_D (y_0)=4t^{1/\alpha} \quad \hbox{ and } \quad
|x-x_0|, \, |y-y_0| \in [t^{1/\alpha},  4t^{1/\alpha}).
$$
By the semigroup property,  with $B:= B(x_0, 4^{-1} t^{1/\alpha})$
and $\wt B:= B(y_0, 4^{-1} t^{1/\alpha})$
\begin{eqnarray*}
p^a_D(t, x, y) &=& \int_{D} p^a_{D}(t/3, x, z) \int_{D}
p^a_{D}(t/3, z, w) p^a_{D}(t/3, w, y) dw dz\\
&\ge& \int_{B} p^a_{D}(t/3, x, z) \int_{\widetilde{B}} p^a_{D}(t/3,
z, w) p^a_{D}(t/3, w, y) dw dz\\
&\ge& \inf_{(z,w) \in B \times \wt {B}} p^a_{D}(t/3, z, w)\int_{B}
p^a_{D}(t/3, x, z) dz\int_{\widetilde{B}} p^a_{D}(t/3, w, y) dw  .
\end{eqnarray*}
Since for $z\in B$ and $w\in \wt B$,
$$\delta_D(z) \ge \delta_D(x_0)-|x_0-z| \ge t^{1/\alpha}, \quad
\delta_D(w) \ge \delta_D(y_0)-|y_0-w| \ge t^{1/\alpha},
$$
$$
|z-w| \le |z-x_0|+|x_0-x|+|x-y|+|y-y_0|+|y_0-w|
<10 t^{1/\alpha}\psi^a(|x-y|)^{1/(d+\alpha)},
$$
by combining Proposition \ref{step1} and Proposition \ref{step3}, we
have that there exists $c_1=c_1(\alpha, \beta, r_0, M)>0$ such that
for all $a\in (0, M]$,
$$
\inf_{(z,w) \in B \times \widetilde{B}} p^a_{D}(t/3, z, w) \,\geq \,
c_1 t^{-d/\alpha}.
$$
Since $\delta_D (x) \leq 2t^{1/\alpha}<r_0/8$ and $\delta_D (y) \leq
3 t^{1/\alpha}$, we have by Lemma \ref{l:st2_3}
$$
p^a_D(t, x, y) \,\ge \, c_2\,
t^{-d/\alpha-1}\,\delta_D(x)^{\alpha/2} \,\delta_D(y)^{\alpha/2}
$$
for some positive constant  $c_2=c_2(\alpha, \beta, M, r_0)$. \qed

\begin{prop}\label{t:lu4} Suppose that $M>0$ and $(t, x, y)\in (0, T_0]
\times D \times D$ with  $\delta_D(x) \le (t/2)^{1/\alpha} \leq
\delta_D(y)$ and $|x-y|^\alpha \ge  t \psi^a
(|x-y|)^{\alpha/(d+\alpha)}$. Then there exists a constant
$C_{29}=C_{29}(\alpha, \beta, M, r_0)>0$ such that for all $a\in (0,
M]$,
\begin{equation}\label{e:lu4}
p^a_D(t, x, y)\,\ge  \,C_{29} \, {t^{1/2} \,
\delta_D(x)^{\alpha/2}}{J^a(x,y)}.
\end{equation}
\end{prop}

\pf Since $\delta_D (x) \leq (t/2)^{1/\alpha} \leq r_0/16$, there
are $z_x\in \partial D$ and $z_0 \in D$ such that
$\delta_D(x)=|x-z_x|$ and $\partial B(z_0, 2t^{1/\alpha}) \cap
\partial D = \{z_x\}$. Choose $x_0$ in
$B(z_0,2t^{1/\alpha})$ and $\kappa=\kappa(\alpha) \in (0, 1)$ such
that
$$
B\left(x_0, 2\kappa t^{1/\alpha}\right)  \subset
B\left(z_0,(2-2^{-2/\alpha}) t^{1/\alpha}\right) \cap
B\left(x,(1-2^{-1-2/\alpha}) t^{1/\alpha}\right).
$$
Such a ball $B(x_0, 2\kappa t^{1/\alpha})$ always exists because
$$
2 \,<\,(2-2^{-1}) + (1-2^{-2})\, <\,(2-2^{-2/\alpha}) +
(1-2^{-1-2/\alpha}).
$$
Since $|x-y| \ge  t^{1/\alpha} \psi^a (|x-y|)^{1/(d+\alpha)}$, we get
that for every  $z \in B(x_0, \kappa t^{1/\alpha})$, $\delta_D(z)\ge
(t/4)^{1/\alpha}$ and
$$
|y-z| \ge |y-x|-|z-x| \ge 2^{-1}
(t/4)^{1/\alpha}\psi^a(|x-y|)^{1/(d+\alpha)}.
$$
On the other hand, for every $z\in B(x_0, \kappa t^{1/\alpha})$,
$$
|z-y| \,\leq\, |z-x|+|x-y| \,\leq\, (1-2^{-1-2/\alpha})
t^{1/\alpha}+|x-y| \,< \,2|x-y|.
$$
Thus by the semigroup property and Propositions \ref{step1} and
\ref{step3}, there exist positive constants $c_1, c_2$ and $c_3$
depending only on $(\alpha, \beta, r_0, M)$ such that for all $a\in
(0,M]$,
\begin{eqnarray*}
p^a_D(t, x, y) &=& \int_{D} p^a_{D}(t/2, x, z) p^a_{D}(t/2, z, y)dz\\
&\ge& \int_{B(x_0, \kappa t^{1/\alpha})} p^a_{D}(t/2, x, z)
p^a_{D}(t/2,
z, y) dz\\
&\ge& c_1 t\int_{B(x_0, \kappa t^{1/\alpha})} p^a_{D}(t/2,
x, z) {J^a(z,y)} dz \\
&\ge &c_2 {t}{j^a(2|x-y|)}\int_{B(x_0, \kappa t^{1/\alpha})}
p^a_{D}(t/2, x, z)dz\\
&\ge& c_3 {t}{j^a(|x-y|)} \, \bP_x \left(X^{a, D}_{t/2} \in B(x_0,
\kappa t^{1/\alpha}) \right).
\end{eqnarray*}
Applying Lemma  \ref{l:st2_3}, we arrive at the conclusion of the
proposition.
 \qed

\begin{prop}\label{t:lu5}
Suppose that $M>0$ and $(t, x, y)\in (0, T_0]\times D\times D$ with
$$
\max \left\{ \delta_D(x), \,  \delta_D(y) \right\} \le
(t/2)^{1/\alpha} \leq |x-y| \psi^a(|x-y|)^{-1/(d+\alpha)}.
$$
Then there exists a constant
$C_{30}=C_{30}(\alpha, \beta, M, r_0)>0$ such that for all $a\in (0,
M]$,
\begin{equation}\label{e:lu5}
p^a_D(t, x, y)\,\ge  \,C_{30}
{\delta_D(x)^{\alpha/2}\delta_D(y)^{\alpha/2}}{J^a(x,y)}.
\end{equation}
\end{prop}

\pf As in the first paragraph of the proof of  Proposition
\ref{step2}, set $z_x\in \partial D$ and $x_0 \in D$ so that
$|x-z_x|=\delta_D (x)$ and $\partial B(x_0, 3t^{1/\alpha}) \cap
\partial D = \{z_x\}$. Let $\kappa :=1-2^{-1/\alpha}$.  Note that
for every $z \in B(x_0, \kappa t^{1/\alpha})$, we have
$$
4t^{1/\alpha}\ge \delta_D(z) \ge 2(t/2)^{1/\alpha}.
$$
If $|y-z|\le t^{1/\alpha}\psi^a(|y-z|)^{1/(d+\alpha)}$, we can apply
Proposition \ref{step2} and the assumption
$$
(t/2)^{1/\alpha} \leq |x-y| \psi^a(|x-y|)^{-1/(d+\alpha)}
$$
to get that
$$
p^a(t/2, z, y)\ge c_1t^{1/2}\, \delta_D(y)^{\alpha/2}
J^a(x,y).
$$
If $|y-z|\ge t^{1/\alpha}\psi^a(|y-z|)^{1/(d+\alpha)}$, we can apply
Proposition \ref{t:lu4} to get that
$$
p^a(t/2, z, y)\ge c_2t^{1/2}\, \delta_D(y)^{\alpha/2}
J^a(y,z).
$$
For $z\in B(x_0, \kappa t^{1/\alpha})$, we have
\begin{eqnarray*}
|z-y|&\le& |x-y|+|x_0-x|+|x_0-z| \,\le\, |x-y| + 4t^{1/\alpha}
\\&\le& |x-y| +
2^{2+1/\alpha}(t/2)^{1/\alpha} \psi^a
(|x-y|)^{1/(d+\alpha)}\\
&\le& (1 + 2^{2+1/\alpha})|x-y|.
\end{eqnarray*}
Thus if $|y-z|\ge t^{1/\alpha}\psi^a(|y-z|)^{1/(d+\alpha)}$, we have
$$
p^a(t/2, z, y)\ge c_3t^{1/2}\, \delta_D(y)^{\alpha/2}
J^a(x,y).
$$
Consequently we have for all $z\in B(x_0, \kappa t^{1/\alpha})$
$$
p^a(t/2, z, y)\ge c_4t^{1/2}\, \delta_D(y)^{\alpha/2}
J^a(x,y).
$$
Hence by the semigroup property we get
\begin{eqnarray*}
p^a_D(t, x, y) &=& \int_{D} p^a_{D}(t/2, x, z)
p^a_{D}(t/2, z, y) dz\\
&\ge& \int_{B(x_0, \kappa t^{1/\alpha})} p^a_{D}(t/2, x, z)
p^a_{D}(t/2, z, y) dz \\
&\ge& c_4 \int_{B(x_0, \kappa t^{1/\alpha})} p^a_{D}(t/2, x,
z) t^{1/2}\, \delta_D(y)^{\alpha/2}J^a(x,y) dz\\
&=&c_4 \, {t^{1/2}\, \delta_D(y)^{\alpha/2}}{J^a(x,y)}\int_{B(x_0,
\kappa t^{1/\alpha})} p^a_{D}(t/2, x, z)dz\\
&=& c_4 \, {t^{1/2} \, \delta_D(y)^{\alpha/2}}{J^a(x,y)}\, \bP_x
\left( X^{a,D}_{t/2}
\in B(x_0, \kappa t^{1/\alpha})  \right)\\
&=& c_5 \, {\delta_D(x)^{\alpha/2}\delta_D(y)^{\alpha/2}}{J^a(x,y)}.
\end{eqnarray*}
We arrive at the conclusion of the proposition. \qed

\noindent {\bf Proof of  Theorem \ref{t:lb}.}
In this proof, for two non-negative functions $f$ and $g$,
the notation $f\asymp g$ means that there are positive constants
$c_1$ and $c_2$ depending only on $M$, $d$, $\alpha$ and $\beta$ so
that $c_1g( x)\leq f (x)\leq c_2 g(x)$ in the common domain of
definition for $f$ and $g$.

We first assume that
$t \le T_0$.
\begin{enumerate}
\item We first consider the case $|x-y| \psi^a
(|x-y|)^{-1/(d+\alpha)}
  \leq
t^{1/\alpha}$.
We claim that in this case
\begin{equation}\label{e:ccc1}
p_D(t, x, y) \ge c t^{-d/\alpha}\left( 1\wedge
\frac{\delta_D(x)^{\alpha/2}}{\sqrt{t}}\right) \left( 1\wedge
\frac{\delta_D(y)^{\alpha/2}}{\sqrt{t}}\right).
\end{equation}
This will be proved by considering the following two possibilities.
  \begin{enumerate}
\item
  $\max\{\delta_D(x), \delta_D(y), |x-y|\psi^a
(|x-y|)^{-1/(d+\alpha)}\} \le t^{1/\alpha}$:
   Proposition \ref{step2} and symmetry yield \eqref{e:ccc1}

\item $\max \{\delta_D(x), \delta_D(y)\}
  \ge t^{1/\alpha} \geq |x-y| \psi^a
(|x-y|)^{-1/(d+\alpha)}$:

If $\max\{\delta_D(x), \delta_D(y)\}\ge t^{1/\alpha} \ge 2|x-y|\psi^a
(|x-y|)^{-1/(d+\alpha)}$,
 \eqref{e:ccc1}  follows from  Proposition \ref{step1}.

If $\min\{\delta_D(x), \delta_D(y)\}\ge t^{1/\alpha}$ and $|x-y|\le
t^{1/\alpha}<2|x-y|\psi^a
(|x-y|)^{-1/(d+\alpha)}$,
$$
 \frac{t\psi^a(|x-y|)}{|x-y|^{d+\alpha}} \asymp  t^{-d/\alpha}\left( 1\wedge
\frac{\delta_D(x)^{\alpha/2}}{\sqrt{t}}\right) \left( 1\wedge
\frac{\delta_D(y)^{\alpha/2}}{\sqrt{t}}\right).
 $$
If
 $\max\{\delta_D(x), \delta_D(y)\}\ge t^{1/\alpha}$,
$\min\{\delta_D(x), \delta_D(y)\} < t^{1/\alpha}$ and $|x-y|\psi^a
(|x-y|)^{-1/(d+\alpha)}\le
t^{1/\alpha}<2|x-y|\psi^a
(|x-y|)^{-1/(d+\alpha)}$,
$$
 \left(
\frac{\delta_D(x)^{\alpha/2}}{\sqrt{t}}\right) \left(
\frac{\delta_D(y)^{\alpha/2}}{\sqrt{t}}\right)\asymp  \left( 1\wedge
\frac{\delta_D(x)^{\alpha/2}}{\sqrt{t}}\right) \left( 1\wedge
\frac{\delta_D(y)^{\alpha/2}}{\sqrt{t}}\right) $$
Thus
by combining Proposition \ref{step3} and Proposition \ref{step2},
we get \eqref{e:ccc1} for
the case of $\max\{\delta_D(x), \delta_D(y)\}\ge t^{1/\alpha}$ and $|x-y|\le
t^{1/\alpha}<2|x-y|\psi^a
(|x-y|)^{-1/(d+\alpha)}$.

\end{enumerate}
\item
Now we consider the case $|x-y|\psi^a
(|x-y|)^{-1/(d+\alpha)} \ge t^{1/\alpha}$ and claim that
\begin{equation}\label{e:ccc2}
p_D(t, x, y) \ge c \left( 1\wedge
\frac{\delta_D(x)^{\alpha/2}}{\sqrt{t}}\right) \left( 1\wedge
\frac{\delta_D(y)^{\alpha/2}}{\sqrt{t}}\right)
 {t}j^a(|x-y|).
\end{equation}
\begin{enumerate}
\item $\min
 \{\delta_D(x), \delta_D(y)\} \le (t/2)^{1/\alpha}$ and $|x-y| \psi^a
(|x-y|)^{-1/(d+\alpha)}\ge t^{1/\alpha}$:
 By symmetry we can assume $\delta_D(x) \le (t/2)^{1/\alpha}$. Thus
combining
Propositions \ref{t:lu4} and \ref{t:lu5},  we have
\eqref{e:ccc2} for this case.
 \item $\min
\{\delta_D(x), \delta_D(y)\} \ge (t/2)^{1/\alpha}$ and $|x-y| \psi^a
(|x-y|)^{-1/(d+\alpha)}\ge
t^{1/\alpha}$.
In this case, clearly
$$ \left( 1\wedge
\frac{\delta_D(x)^{\alpha/2}}{\sqrt{t}}\right) \left( 1\wedge
\frac{\delta_D(y)^{\alpha/2}}{\sqrt{t}}\right) \asymp
\left(
\frac{\delta_D(x)^{\alpha/2}}{\sqrt{t}}\right) \left(
\frac{\delta_D(y)^{\alpha/2}}{\sqrt{t}}\right).
 $$
 Thus
Proposition \ref{step3}
 yields \eqref{e:ccc2}.
\end{enumerate}
We have arrived at the conclusion of Theorem
\ref{t:lb} for $t\le T_0$.
\end{enumerate}

Assume $T=2T_0$. Recall that $T_0=(r_0/16)^\alpha$. For $(t, x,
y)\in (T_0, 2T_0] \times D \times D$, let $x_0, y_0\in D$ be such
that $\max\{|x-x_0|, |y-y_0|\} <r_0$ and $\min\{\delta_D (x_0),
\delta_D (y_0)\}\geq r_0/2$. Note that, since
\begin{equation}\label{e:jjj}
j^a(r)\le c_1 j^a(2r), \quad \text{ for all } r>0,
\end{equation}
if $|x-y| \ge 4r_0$, then
$\frac12 |x-y| \le |x-y|-2r_0 \le |x_0-y_0| \le |x-y|+2r_0 \le \frac32 |x-y|$, and so
$c_2^{-1} J^a(x_0, y_0) \le J^a(x, y) \le c_2 J^a(x_0, y_0)$ for
some constant $c_2=c_2(M)>1$. Thus by
 considering the cases $|x-y| \ge 4r_0$ and $|x-y| < 4r_0$, we have
\begin{equation}\label{e:N1}
(t/2)^{-d/\alpha} \wedge \frac{tJ^a(x_0, y_0)}2 \geq c_3 \left(
t^{-d/\alpha} \wedge (tJ^a(x, y))\right) .
\end{equation}
Similarly, there is a positive constant $c_2$ such that
\begin{eqnarray}
(t/3)^{-d/\alpha} \wedge \frac{tJ^a(x, z)}3
&\ge& c_4 \left( (t/(12))^{-d/\alpha} \wedge
\frac{t J^a(x_0,z)}{12} \right), \quad z\in D, \nn \\
 (t/3)^{-d/\alpha} \wedge \frac{tJ^a(w,y)}{3}
&\ge& c_4 \left( (t/(12))^{-d/\alpha} \wedge
\frac{t J^a(w,y_0)}{12} \right), \quad w\in D. \label{e:N2}
\end{eqnarray}
By \eqref{e:N2} and the lower bound estimate in Theorem \ref{t:lb}
for $p^a_D$ on $(0, T_0]\times D\times D$, we have
\begin{align*}
&p^a_D(t, x, y) \,= \, \int_{D\times D}
p^a_D(t/3, x, z)p^a_D (t/3, z, w) p^a_D(t/3, w, y) dz dw \\
&\geq \, c_5 \left(1\wedge
\frac{\delta_D(x)^{\alpha/2}}{\sqrt{t/3}}\right) \left( 1\wedge
\frac{\delta_D(y)^{\alpha/2}}{\sqrt{t/3}}\right) \int_{D\times D}
\left( (t/3)^{-d/\alpha} \wedge \frac{tJ^a(x, z)}3 \right)
\left(1\wedge
\frac{\delta_D(z)^{\alpha/2}}{\sqrt{t/3}}\right)  \\
& \hskip 0.6truein \cdot  p^a_D(t/3, z, w) \left( (t/3)^{-d/\alpha}
\wedge \frac{t J^a(w,y)}3 \right) \left(1\wedge
\frac{\delta_D(w)^{\alpha/2}}{\sqrt{t/3}}\right) dzdw \\
&\geq \, c_6 \left(1\wedge
\frac{\delta_D(x)^{\alpha/2}}{\sqrt{t}}\right) \left( 1\wedge
\frac{\delta_D(y)^{\alpha/2}}{\sqrt{t}}\right) \int_{D\times D}
\left( \left(\frac{t}{12}\right)^{-d/\alpha} \wedge
\frac{tJ^a(x_0,z)}{12} \right)\left(1\wedge
\frac{\delta_D(z)^{\alpha/2}}{\sqrt{t/3}}\right)  \\
& \hskip 0.6truein \cdot \,  p^a_{D}(t/3, z, w) \left(
\left(\frac{t}{12}\right)^{-d/\alpha} \wedge \frac{tJ^a(w,y_0)}{12}
\right)\left(1\wedge
\frac{\delta_D(w)^{\alpha/2}}{\sqrt{t/3}}\right) dzdw
\end{align*}
for some positive constants $c_i, i=3, 4$.
Let $D_1:=\{z \in D: \delta_D(z) > r_0/4\}$. Clearly, $x_0, y_0\in
D_1$ and
\begin{equation}\label{e:nem}
\min\{\delta_{D_1} (x_0), \delta_{D_1} (y_0)\}\geq
r_0/4=4(T_0)^{1/\alpha} \ge 4 (t/2)^{1/\alpha}.
\end{equation}
By \eqref{e:1.4} and \eqref{e:N1}, we have
\begin{align*}
&\int_{D\times D}
\left(\left(\frac{t}{12}\right)^{-d/\alpha} \wedge
\frac{t  J^a(x_0,z)}{12} \right)\left(1\wedge
\frac{\delta_D(z)^{\alpha/2}}{\sqrt{t/3}}\right)  \\
& \hskip 0.6truein \cdot \,  p^a_{D}(t/3, z, w)
  \left(\left(\frac{t}{12}\right)^{-d/\alpha} \wedge
\frac{t J^a(w,y_0)}{12} \right)\left(1\wedge
\frac{\delta_D(w)^{\alpha/2}}{\sqrt{t/3}}\right) dzdw\\
 \ge & c_7\int_{D_1\times D_1}
\left( \left(\frac{t}{12}\right)^{-d/\alpha} \wedge
\frac{t  J^a(x_0,z)}{12} \right)  p^a_{D}(t/3, z, w)
  \left( \left(\frac{t}{12}\right)^{-d/\alpha} \wedge
\frac{t J^a(w,y_0)}{12} \right) dzdw\\
\ge& \,c_8\,
 \int_{D_1\times D_1}
 p^a (t/(12), x_0, z) p^a_{D_1} (t/3, z, w)p^a (t/(12), w, y_0) dzdw \\
\geq &\,  c_8 \, \int_{D_1 \times D_1 }
  p^a_{D_1} (t/(12), x_0, z) p^a_{D_1} (t/3, z, w)p^a_{D_1} (t/(12), w, y_0) dzdw \\
=&\, c_8\,
p^a_{D_1}(t/2, x_0, y_0) \geq  c_9
  \left(  (t/2)^{-d/\alpha}
\wedge \frac{t J^a(x_0, y_0)}2 \right) \geq  c_{10}
  \left( t^{-d/\alpha} \wedge (tJ^a(x, y)) \right)\end{align*}
for some positive constants $c_i, i=7, \cdots, 10$. Here both Propositions
\ref{step1} and \ref{step3} are used in the third inequality in view of
\eqref{e:nem}. By repeating the argument above, we have proved
Theorem \ref{t:lb}.\qed

\noindent {\bf Proof of  Theorem \ref{t:main}.} Theorems \ref{t:ub}
and \ref{t:lb} give Theorem \ref{t:main}(i). By \cite{Gr, KS}, for
any bounded open set $D$ in $\R^d$, $X^{a, D}$ is intrinsically
ultracontractive. Since the function $\psi^a(|x-y|)$ is bounded
above and below by a positive constant if $D$ is bounded, using the
intrinsic ultracontractivity of $X^a$ on bounded open set and the
continuity of eigenvalues proved in \cite{CS8}, the proof
 of Theorem \ref{t:main} (ii) is almost identical to the one of
\cite[Theorem 1.1(ii)]{CKS2}. We omit the details. \qed

\noindent {\bf Proof of  Corollary \ref{C:1.2}.}
The lower bound estimate in \eqref{e:G} follows from
\eqref{gfcrln} and Theorem \ref{T:gfcesks}.

Since the function $\psi^a(|x-y|)$ is bounded above and below by a
positive constant if $D$ is bounded, by integrating the two-sided
heat kernel estimates in Theorem \ref{t:main} with respect to $t$,
the proof of the upper bound estimate in \eqref{e:G} is identical to
the one of \cite[Corollary 1.2]{CKS} so we omit its details here.
\qed

\begin{thm}[Uniform boundary Harnack principle]\label{ubhp}
Suppose $M, R \in (0, \infty)$ and that $D$ is an open set in
$\bR^d$, $z\in \partial D$, $r\in (0, R)$ and that $B(A, \kappa r)
\subset D\cap B(z, r)$. There exists
$C_{31}=C_{31}(d, \alpha, \beta, \kappa, M, R)>1$ such that for
every $a \in (0, M]$, and any functions $u, v\ge 0$ on $\bR^d$,
positive regular harmonic for $X^{a}$ in $D\cap B(z, 2r)$ and
vanishing on $D^c\cap B(z, 2r)$, we have
$$
C^{-1}_{31}\frac{u(A)}{v(A)}\le \frac{u(x)}{v(x)} \le C_{31}
\frac{u(A)}{v(A)}, \qquad x\in D\cap B(z,  r).
$$
\end{thm}
\pf Note that by the approximate scaling property in
\eqref{e:scaling}, we have for every $r>0$.
\begin{equation}\label{e:scG1}
G^a_{B(0,r)}(x,y)=r^{\alpha-d} G_{B(0,1)}^{ar^{(\alpha-\beta)/\beta}}(x/r,y/r).
\end{equation}
Thus applying \cite[Corollary 1.2]{CKS} and our Corollary \ref{C:1.2} to \eqref{e:scG1}, we have that for every $R, M>0$, there exists
$c=c(\alpha, \beta, R, M)>0$ such that, for every $a \in (0, M]$ and $0<r \le R$
\begin{equation}\label{e:scG2}
c^{-1} G_{B(x_0,r)}(x,y) \le G^a_{B(x_0,r)}(x,y) \le c G_{B(x_0,r)}(x,y), \quad \forall x,y \in B(x_0,r).
\end{equation}
Using \eqref{e:scG2}, we can get uniform estimates on the Poisson
kernel
$$
K_{B(x_0, r)}^a (x,z):= \int_{  B(x_0, r)} G^a_{B(x_0, r)}(x,y)J^a(y,z)dy
$$
of $B(x_0,r)$ with respect to $X^a$ for $r\in (0, R]$. In particular,  for $ r <|z-x_0|< 2 R$, $K_{B(x_0, r)}^a (x,z)$ is comparable to $K_{B(x_0,
r)}(x,z)$,  the Poisson kernel of $B(x_0,r)$ with respect to $X$
for $r\in (0, R]$. Then using the uniform estimates
on $K_{B(x_0, r)}^a (x,z)$ and \eqref{e:scG2} we
 can easily see that \cite[Lemma 3.3]{SW} can be proved in the same
way. Using the uniform estimates on the Poisson kernel of
$B(x_0,r)$,  \eqref{e:jjj} and \eqref{e:scG2}
we can adapt the argument in
\cite{B, KSV, SW}
to get our
uniform boundary Harnack principle. We omit the details.
\qed

\noindent {\bf Proof of  Theorem \ref{t:bhp}}
First we observe that Harnack inequality holds for the process $X:=X^1$ by \cite{RSV}.
That is,
there exists a constant $c_1=c_1(\alpha, \beta, M)>0$
such that for any $r\in (0, M^{\beta/(\alpha-\beta)}]$, $x_0\in \R^d$ and
any function $v \ge 0$ harmonic in
$B(x_0, r)$ with respect to
 $X$,
 we have
\begin{equation}\label{e:HP1}
v(x)\le c_1v(y) \qquad \mbox{ for all } x, y\in B(x_0, r/2).
\end{equation}
Note that  for any $a\in (0, M]$, $X^a$ has the same
distribution as $\{\lambda X_{\lambda^{-\alpha}t}, t\geq 0\}$,
where $\lambda = a^{\beta/(\beta-\alpha)} \geq M^{\beta/(\beta-\alpha)}$.
Consequently, if $u$ is harmonic in $B(x_0, r)$ with respect to  $X^a$,
where $r \in (0,1]$, then $v(x):=u(\lambda x)$ is harmonic in $B(\lambda^{-1}x_0,
\lambda^{-1}r)$ with respect to $X$ and $\lambda^{-1}r \le
M^{\beta/(\beta-\alpha)} $. So by \eqref{e:HP1}
$$
u(\lambda x)=v(x)\le c_1v(y)= c_1u(\lambda y)
\qquad \mbox{ for all } x, y\in B(\lambda^{-1}x_0, \lambda^{-1} r/2).
$$
That is,
\begin{equation}\label{e:HP2}
u( x)\le c_1u( y) \qquad \mbox{ for all } x, y\in B(x_0, r/2).
\end{equation}
In other words, uniform Harnack inequality holds (for every $r \le
1$) for the family of processes $\{X^a, a\in (0, M]\}$.

Since $D$ is $C^{1,1}$ open set, there exists $r_0 \le  R_0$ such
that the following holds: for every $Q\in \partial D$ and $r\le r_0$
there is a ball $B=B(z^r_Q, r)$ of radius $r$ such that $B\subset
\bR^d \setminus \overline D$ and $\partial B \cap \partial D=\{Q\}$.
In addition, it follows \cite[Lemma 2.2]{S} that, for each $Q \in
\partial D$, we can choose a constant $c_2=c_2(d, \Lambda) \in
(0,1/8]$ and a bounded $C^{1,1}$ open set $U_Q$ with uniform
characteristics $(R_*, \Lambda_*)$ depending on $(R_0, \Lambda)$
such that $B(Q, c_2 r_0) \cap D \subset U_Q \subset B(Q, r_0) \cap
D$ and
\begin{equation}\label{e:dcp}
\delta_D(y)=\delta_{U_Q}(y) \quad \text{ for every }y \in B(Q, c_2 r_0) \cap D.
\end{equation}

Assume $a \in [0, M]$, $r \in (0, c_2r_0]$, $Q\in \partial D$ and $u$ is
nonnegative function in $\R^d$ harmonic in $D \cap B(Q,
r)$ with respect to $X^{a}$ and vanishes continuously on $ D^c \cap
B(Q, r)$. Let $z_Q:=z_Q^{c_2r_0}$.
By the boundary Harnack principle (Theorem \ref{ubhp}), there exists
a constant $c_3=c_3(\alpha, \beta, a, R_0, \Lambda, M)$ such that
$$
\frac{u(x)}{u(y)} \le c_3 \frac{G^a_{U_Q}(x, z_Q)}{G^a_{U_Q}(x,
z_Q)} \quad \text{ for every } x,y \in B(Q, r/8) \cap D.
$$
Now applying Corollary \ref{C:1.2} to $G^a_{U_Q}(x, z_Q)$ and
$G^a_{U_Q}(x, z_Q)$, then using \eqref{e:dcp}, we conclude that
\begin{equation}\label{e:bbhp}
\frac{u(x)}{u(y)} \le c_4
\frac{\delta^{\alpha/2}_{U_Q}(x)}{\delta^{\alpha/2}_{U_Q}(y)} = c_4
\frac{\delta^{\alpha/2}_D(x)}{\delta^{\alpha/2}_D(y)}\quad \text{
for every } x,y \in B(Q, c_2 r) \cap D
\end{equation}
for some $c_4=c_4(\alpha, \beta, a, R_0, \Lambda, M)>0$.

Now Theorem \ref{t:bhp} follows from
the uniform Harnack principle in \eqref{e:HP2}, \eqref{e:bbhp} and a standard chain
argument.
\qed

\vskip 0.3truein

{\bf Zhen-Qing Chen}

Department of Mathematics, University of Washington, Seattle,
WA 98195, USA

E-mail: \texttt{zchen@math.washington.edu}

\bigskip

{\bf Panki Kim}

Department of Mathematical Sciences,
Seoul National University,
San56-1 Shinrim-dong Kwanak-gu,
Seoul 151-747, Republic of Korea

E-mail: \texttt{pkim@snu.ac.kr}

\bigskip

{\bf Renming Song}

Department of Mathematics, University of Illinois, Urbana, IL 61801, USA

E-mail: \texttt{rsong@math.uiuc.edu}


\bigskip

\begin{thebibliography}{99}

\bibitem{B} K. Bogdan, The boundary Harnack principle for the fractional
Laplacian. {\it Studia Math. \bf 123} (1997), 43--80.

\bibitem{BBC} K. Bogdan, K. Burdzy and Z.-Q. Chen,
Censored stable processes, {\em Probab. Theory Related Fields}, {\bf
127} (2003),  89--152.

\bibitem{CKS}  Z.-Q. Chen, P. Kim, and R. Song,
Heat kernel estimates for Dirichlet fractional Laplacian. {\em J.
European Math. Soc.}, (to appear), 2009.

\bibitem{CKS1} Z.-Q. Chen, P. Kim and R. Song,
Two-sided heat kernel estimates for censored stable-like processes.
{\em Probab. Theory Relat. Fields}, (to appear), 2009.

\bibitem{CKS2} Z.-Q. Chen, P. Kim and R. Song,
Sharp heat kernel estimates for relativistic stable processes in open sets.
Preprint 2009,
arXiv:0908.1509 [math.PR].

\bibitem{CKSV} Z.-Q. Chen, P. Kim, R. Song and Z. Vondra\v cek,
Boundary Harnack pinciple for $\Delta + \Delta^{\alpha/2}$.
Preprint, 2009,
arXiv:0908.1559 [math.PR].

\bibitem{CK} Z.-Q. Chen and T. Kumagai,
Heat kernel estimates for stable-like processes on $d$-sets, {\em
Stoch. Proc. Appl.},  {\bf 108} (2003), 27--62.

\bibitem{CK2}  Z.-Q. Chen and  T. Kumagai,
Heat kernel estimates for jump processes of mixed types on metric
measure spaces. {\it Probab. Theory Relat. Fields},  {\bf 140}
(2008), 277--317.



\bibitem{CS1} Z.-Q.~Chen and R.~Song,
\newblock Estimates on Green functions and Poisson kernels of
symmetric stable processes,
\newblock {\em Math. Ann.}, {\bf 312} (1998), 465-601.

\bibitem{CS8} Z.-Q. Chen and R. Song, Continuity of eigenvalues for
subordinate processes in domains. {\it Math. Z.}, {\bf 252} (2006),
71--89.

\bibitem{D1}
E. B. Davies, Explicit constants for Gaussian upper bounds on heat
kernels. {\em Amer. J. Math.} {\bf 109} (1987), 319--333.

\bibitem{D2}
E. B. Davies, The equivalence of certain heat kernel and Green
function bounds. {\em J. Funct. Anal.} {\bf 71} (1987), 88--103.

\bibitem{Deb}  E. B. Davies,
{\em Heat Kernels and Spectral Theory}, Cambridge University Press,
Cambridge, 1989.

\bibitem{DS} E. B. Davies and B. Simon, Ultracontractivity and heat kernels
for Schr\"odinger operator and Dirichlet Laplacians. {\it J. Funct. Anal.
\bf 59} (1984), 335-395.


\bibitem{Gr}
T. Grzywny, Intrinsic ultracontractivity for L\'evy processes. {\em
Probab. Math. Statist.} {\bf 28(1)} (2008), 91--106.

\bibitem{GR} T. Grzywny and M. Ryznar,
Estimates of Green functions for some perturbations of fractional
Laplacian. {\em Illinois J. Math.} {\bf 51}  (2007), 1409--1438.

\bibitem{G}  Q.-Y. Guan, Boundary Harnack inequality for regional
fractional Laplacian.  arXiv:0705.1614v2 [math.PR]



\bibitem{KS} P. Kim and  R. Song,
Intrinsic Ultracontractivity for Non-symmetric L\'evy Processes.
{\em Forum Math.} {\bf 21(1)} (2009), 43--66.
Erratum to: Intrinsic ultracontractivity for non-symmetric L\'evy
processesÓ [Forum Math. 21 (2009) 43--66], to appear in {\em Forum
Math.}

\bibitem{KSV} P. Kim, R. Song and Z. Vondra\v{c}ek, Boundary Harnack
principle for subordinate Brownian motion. {\em Stoch. Proc. Appl.},
{\bf 119} (2009), 1601--1631.



\bibitem{Ku1} T. Kulczycki,
Properties of Green function of symmetric stable processes, {\em
Probab. Math. Stat.} {\bf 17}(1997), 381--406.


\bibitem{RSV}
M. Rao, R. Song and Z. Vondra\v{c}ek,
Green function estimates and Harnack inequality for
subordinate Brownian motions.
 {\em Potential Anal.} {\bf 25} (2006), 1--27

\bibitem{R} M. Ryznar, Estimates of Green function for relativistic
$\alpha$-stable process.  {\em Potential Anal. } {\bf 17} (2002),
1--23.

\bibitem{S} R. Song, Estimates on the Dirichlet heat kernel of
domains above the graphs of bounded $C\sp {1,1}$ functions. {\em
Glas. Mat.}  {\bf 39} (2004), 273--286.


\bibitem{SV} R. Song and Z. Vondracek, Sharp bounds for Green functions
and jumping functions of subordinate killed Brownian motions in
bounded $C\sp {1,1}$ domains.  {\em Electron. Comm. Probab.}, {\bf
9} (2004), 96--105.

\bibitem{SV06} R.~Song, Z.~Vondracek, Potential theory of special subordinators
and subordinate killed stable processes.  {\em J. Theoret. Probab.},
{\bf 19} (2006), 817--847.

\bibitem{SV08} R.~Song, Z.~Vondracek,  On the relationship between
subordinate killed and killed subordinate processes, \emph{Elect.
Commun.  Probab.}, {\bf 13} (325--336) 2008.


\bibitem{SW}
R. Song and J. Wu, Boundary Harnack principle for symmetric stable
processes, {\it J. Funct. Anal.} {168(2)} (1999), 403--427.

\bibitem{Sz1} P. Sztonyk,
On harmonic measure for L\'evy processes,
{\em Probab. Math. Statist.}, {\bf 20} (2000),
383--390.


\bibitem{Zq3} Q. S. Zhang, The boundary behavior of heat kernels of Dirichlet
Laplacians, {\em J. Differential Equations}, {\bf 182} (2002),
416--430.

\end{thebibliography}
\end{document}